
\def\forces{\parallel\!\!\! -}


\def\hexnumber#1{\ifcase#1 0\or1\or2\or3\or4\or5\or6\or7\or8\or9\or
	A\or B\or C\or D\or E\or F\fi }

\font\teneuf=eufm10
\font\seveneuf=eufm7
\font\fiveeuf=eufm5
\newfam\euffam
\textfont\euffam=\teneuf
\scriptfont\euffam=\seveneuf
\scriptscriptfont\euffam=\fiveeuf
\def\frak{\fam\euffam \teneuf}

\font\tenmsx=msam10
\font\sevenmsx=msam7
\font\fivemsx=msam5
\font\tenmsy=msbm10
\font\sevenmsy=msbm7
\font\fivemsy=msbm5
\newfam\msxfam
\newfam\msyfam
\textfont\msxfam=\tenmsx  \scriptfont\msxfam=\sevenmsx
  \scriptscriptfont\msxfam=\fivemsx
\textfont\msyfam=\tenmsy  \scriptfont\msyfam=\sevenmsy
  \scriptscriptfont\msyfam=\fivemsy
\edef\msx{\hexnumber\msxfam}

\mathchardef\upharpoonright="0\msx16
\let\restriction=\upharpoonright
\def\Bbb#1{\tenmsy\fam\msyfam#1}

\def\re{{\restriction}}

\def\Smallskip{\vskip1.2truecm}

\def\Bigskip{\vskip2.2truecm}
\def\Hoskip{\hskip0.8truecm}

\def\qed{{\vcenter{\hrule height.4pt \hbox{\vrule width.4pt height5pt
 \kern5pt \vrule width.4pt} \hrule height.4pt}}}
\def\ok{\vbox{\hrule height 8pt width 8pt depth -7.4pt
    \hbox{\vrule width 0.6pt height 7.4pt \kern 7.4pt \vrule width 0.6pt height 7.4pt}
    \hrule height 0.6pt width 8pt}}
\def\nt{{\leq}\kern-1.5pt \vrule height 6.5pt width.8pt depth-0.5pt \kern 1pt}
\def\sd{{\times}\kern-2pt \vrule height 5pt width.6pt depth0pt \kern1pt}
\def\zp#1{{\hochss Y}\kern-3pt$_{#1}$\kern-1pt}

\def\bb{{\frak b}}
\def\cc{{\frak c}}
\def\dd{{\frak d}}

\def\hh{{\frak h}}
\def\hhom{{\frak{hom}}}

\def\pp{{\frak p}}
\def\ppar{{\frak{par}}}

\def\FFF{{\frak F}}
\def\rr{{\frak r}}
\def\ss{{\frak s}}

\def\uu{{\frak u}}

\def\CC{{\Bbb C}}

\def\LL{{\Bbb L}}

\def\MM{{\Bbb M}}

\def\PP{{\Bbb P}}
\def\QQ{{\Bbb Q}}

\def\A{{\cal A}}

\def\D{{\cal D}}

\def\F{{\cal F}}
\def\G{{\cal G}}

\def\I{{\cal I}}

\def\P{{\cal P}}

\def\T{{\cal T}}
\def\U{{\cal U}}
\def\V{{\cal V}}
\def\W{{\cal W}}

\def\add{{\sanse add}}
\def\cov{{\sanse cov}}

\def\non{{\sanse non}}
\def\cof{{\sanse cof}}
\def\Spec#1{{\sanse Spec}$(#1)$}
\def\Specs#1{{\sanse Spec}$^\star(#1)$}
\def\Add#1{{\sanse add}$({\cal #1})$}
\def\Cov#1{{\sanse cov}$({\cal #1})$}
\def\Non#1{{\sanse non}$({\cal #1})$}

\def\Cof#1{{\sanse cof}$({\cal #1})$}
\def\Cofu#1{{\sanse cof}$(\omom/{\cal #1})$}

\def\sm{{\smallskip}}
\def\ce#1{{\centerline{#1}}}
\def\no{{\noindent}}
\def\la{{\langle}}
\def\ra{{\rangle}}
\def\sub{\subseteq}

\def\ha{{\hat{\;}}}
\def\em{{\emptyset}}
\def\sem{\setminus}
\def\omom{{\omega^\omega}}
\def\omlom{{\omega^{<\omega}}}
\def\omup{{\omega^{\uparrow\omega}}}
\def\omlup{{\omega^{\uparrow < \omega}}}
\def\omoms{{[\omega]^\omega}}
\def\omloms{{[\omega]^{<\omega}}}

\def\Loleriar{\Longleftrightarrow}
\font\small=cmr8 scaled\magstep0
\font\smalli=cmti8 scaled\magstep0
\font\capit=cmcsc10 scaled\magstep0
\font\capitg=cmcsc10 scaled\magstep1

\font\dunhg=cmdunh10 scaled\magstep1
\font\dunhgg=cmdunh10 scaled\magstep2

\font\sanse=cmss10 scaled\magstep0

\font\bolds=cmssdc10 scaled\magstep0

\overfullrule=0pt
\openup1.5\jot

\ce{}
\bigskip
\ce{\dunhgg Ultrafilters on $\omega$ ---}
\bigskip
\ce{\dunhgg --- their ideals and their cardinal characteristics   }
\footnote{}{{\openup-6pt {\small {\smalli
1991 Mathematics subject classification.}
03E05 03E35     \par
{\smalli Key words and phrases.} ultrafilter, P--point,
Ramsey ultrafilter,
character, $\pi$--character, ideal, Ramsey null sets (nowhere Ramsey
sets), cardinal coefficient,
Mathias forcing, Laver forcing, Easton forcing.
\endgraf}}}
\Bigskip
\ce{\capitg J\"org Brendle$^{1,}$\footnote{$^\star$}
{{\small Research partially supported by DFG--grant Nr. Br 1420/1--1.}}
and Saharon Shelah$^{2,}$\footnote{$^{\star\star}$}
{{\openup-6pt {\small Supported by the
   German--Israeli Foundation for Scientific Research \& Development
   Grant No. G-294.081.06/93. \endgraf}}}$^,$\footnote{$^\dagger$}
{{\small Publication 642}}}
\Smallskip
\no {}$^1$  Department of Mathematics, Dartmouth College, 
Bradley Hall, Hanover NH 03755, USA;
Graduate School of Science and Technology, Kobe University,
Rokko--dai, Nada, Kobe 657, Japan;
email: {\sanse jobr@michelangelo.mathematik.uni-tuebingen.de}
\bigskip
\no {}$^2$ Institute of Mathematics, The Hebrew University, Jerusalem 
91904, Israel; Department of Mathematics, Rutgers University,
New Brunswick NJ 08903, USA

\Bigskip
\ce{\capit Abstract}
\bigskip
\no For a free ultrafilter $\U$ on $\omega$ we study
several cardinal characteristics which describe part
of the combinatorial structure of $\U$. We  provide various consistency
results; e.g. we show how to force simultaneously many
characters and many $\pi$--characters. We also investigate
two ideals on the Baire space $\omom$ naturally related
to $\U$ and calculate cardinal coefficients of these
ideals in terms of cardinal characteristics of the
underlying ultrafilter.  

\vfill\eject

{\dunhg Introduction}
\Smallskip

\no Let $\U$ be a non--principal ultrafilter on the natural numbers
$\omega$. Recall that $\U$ is a {\it $P$--point} iff for all
countable $\A\sub\U$ there is $U\in\U$ with $U \sem A$ being finite for all 
$A\in\A$. $\U$ is said to be {\it rapid} iff for all $f\in\omom$
there is $U\in\U$ with $|U \cap f(n) | \leq n$ for all $n\in\omega$.
$\U$ is called {\it Ramsey} iff given any partition $\la A_n ;
\; n\in\omega\ra$ of $\omega$, there is either $n\in\omega$ with
$A_n \in\U$ or $U\in\U$ with $|A_n \cap U | \leq 1$ for all
$n\in\omega$. It is well--known (and easily seen) that Ramsey ultrafilters are
both rapid and $P$--point.

With $\U$ we can associate ideals on the real numbers (more exactly, on the
Baire space $\omom$) in various ways. One way of doing this results
in the well--known ideal $r_\U^0$ of {\it Ramsey null sets}
with respect to $\U$
(see $\S$ 2 for the definition). Another, less known, ideal related
to $\U$ was introduced by Louveau in [Lo] and shown to coincide with
both the meager and the nowhere dense ideals on $\omom$ with
respect to a topology somewhat finer than the standard topology
(see $\S$ 3 for details). This ideal which we call $\ell^0_\U$ is
related to {\it Laver forcing with} $\U$, $\LL_\U$ [Bl 1], in a way
similar to the connection between $r^0_\U$ and {\it Mathias forcing with} $\U$,
$\MM_\U$. Furthermore, $\ell^0_\U$ and $r^0_\U$ coincide in case
$\U$ is a Ramsey ultrafilter [Lo], as do $\LL_\U$ and $\MM_\U$
[Bl 1].

A natural problem which has, in fact, been studied for many 
ideals $\I$ on the reals [BJ 1] is to figure out the relationship between
certain {\it cardinal coefficients} of $\I$ as well as to determine
their possible
values. An example of such a cardinal coefficient is the {\it
additivity} of $\I$, \Add I, that is, the size of the smallest
subfamily of $\I$ whose union is not in $\I$; another one,
the {\it uniformity} of $\I$, \Non I, is the cardinality
of the least set of reals which does not belong to $\I$ (see $\S$ 2 for 
more such coefficients). One of the goals of this work is to carry
out such an investigation for $\I$ being either $\ell^0_\U$
or $r^0_\U$. (In fact, this was the original motivation for this
paper.)

In sections 2 and 3 of the present paper 
we reduce this problem to a corresponding problem about
cardinal characteristics of the underlying ultrafilter $\U$, by actually 
calculating the ideal coefficients in terms of the latter as well as of two
other cardinal invariants of the continuum, the {\it unbounding
number} $\bb$ and the {\it dominating number} $\dd$ (see $\S$ 1
for the definitions). Here, by 
a {\it cardinal characteristic} of $\U$, we mean a cardinal number
describing part of the combinatorial structure of $\U$,
such as  the {\it character}
of $\U$, $\chi (\U)$, that is, the size of the smallest subfamily
$\F$ of $\U$ such that each member of $\U$ contains some
member of $\F$ --- or the {\it $\pi$--character} of $\U$, $\pi
\chi (\U)$, the cardinality of the least $\F \sub \omoms$
such that each element of $\U$ contains an element of $\F$ 
(see $\S$ 1 for details).  We show for example that
\non$(r^0_\U) = \pi\chi (\U)$ (Theorem 1(c) in $\S$ 2)
or that the uniformity of $\ell^0_\U$ can be expressed as the
maximum of $\dd$ and some cardinal closely related to
$\pi\chi (\U)$ (Theorem 2(c) in $\S$ 3). The interest of such
characterizations lies in the fact that, unlike the ideal coefficients,
the ultrafilter characteristics have been studied previously,
in particular in connection with ongoing research on $\beta
\omega$ (see e.g. [vM]) but also in investigations of
the cofinality of ultraproducts of the form $\omom / \U$, 
and so already established results on the latter
can be used to show something on the former. Furthermore,
the ultrafilter characteristics as well as the classical cardinal
invariants of the continuum are combinatorially simpler objects
than the ideal coefficients and thus easier to calculate in any given model
of set theory. Accordingly, we investigate the ultrafilter
characteristics in the remainder of our work ($\S\S$ 1, 4 -- 7).

It turns out that only rather elementary facts about these
characteristics and their relationship to other cardinal invariants
can be proved in $ZFC$. Most of these results which we expound
in section 1 are well--known. To make our paper self--contained, we include
proofs. (For the consequences of these $ZFC$--results on the
ideal coefficients, see the corollaries in sections 2 and 3.)

This leaves the field wide open for independence results of various sorts
to which the main body of the present paper (sections 4 to 7) is
devoted. --- First, we deal with distinguishing between different 
coefficients for a fixed Ramsey ultrafilter $\U$. Most questions one
would ask in this direction have been solved long ago (see $\S\S$ 1 and 4).
The remaining case, to force a Ramsey ultrafilter $\U$ with
$\pi\chi (\U) < \chi (\U)$, is taken care of in a rather straightforward
construction in Theorem 3 in section 4. --- Next, we are concerned with
producing simultaneously many different ultrafilters for which a fixed
cardinal characteristic assumes many different values. For one of our
cardinals, this has been done by Louveau ([Lo], see also $\S$ 1)
under $MA$ long ago. For the others, it is a much more difficult problem
which we tackle in sections 5 and 6. For example we show that given a
set of uncountable cardinals $R$ in a model of $GCH$, we can force
that for each $\lambda \in R$ there is an ultrafilter $\U$ with
$\pi\chi (\U) = \lambda$ (Theorem 4(a) and Corollary 5.5). Similarly,
given a set of cardinals of uncountable cofinality $R$, such a model
can be extended to one which has an ultrafilter (even a 
$P$--point) $\U$ with $\chi (\U)
=\lambda$ for all $\lambda \in R$ (Theorem 5 and Corollary 6.1).
For quite many years, R. Frankiewicz, S. Shelah and P. Zbierski
have planned to write a paper proving this for regulars (i.e. for
any set of regulars $R$, there is a forcing extension with a
$P$--point with character $\lambda$ for each $\lambda \in R$).
The proof of Theorem 5
can be extended in various ways, e.g.  to make all the ultrafilters
Ramsey (Corollary 6.2) or  to prove a dual
result (Theorem 7). It is an elegant
combination of a $ccc$--iteration and an Easton product. Results on
characters and $\pi$--characters like those described in sections
5 and 6 are interesting not just because they shed light on the
ideal coefficients studied in section 2 and 3, but also because
$\chi$ and $\pi\chi$ play a role in the topological investigation of 
$\beta\omega$ (see [vM]). --- Finally, we explore in section 7 the
connection between the ultrafilter characteristics and the {\it reaping}
and {\it splitting numbers} $\rr$ and $\ss$ (see $\S$ 1 for the
definitions). Using iterated forcing we show (Theorem 8) that a result of 
Balcar and Simon ([BS], see also Proposition 7.1) which says that
$\rr$ is the minimum of the $\pi$--characters cannot be dualized
to a corresponding statement about $\ss$. The main
technical device of the proof is a careful analysis of
$\LL_\U$--names for reals where $\U$ is a Ramsey ultrafilter.

We close with a list of open problems in section 8.

All sections of this work from section 2 onwards 
depend on section 1, but can be
read independently of each other; however, $\S$ 3 uses the
basic definitions of $\S$ 2; and sections 5 and 6 are closely 
intertwined.
\bigskip

{\bolds Notational remarks and some prerequisites.} 
We refer to standard texts like [Je]
or [Ku] for any undefined notion. $\cc$ stands for the
cardinality of the continuum. $cf (\kappa)$ is the cofinality
of the cardinal $\kappa$. Given a function $f$, $dom (f)$ is its
domain, $rng (f)$ its range, and if $A \sub dom (f)$, then $f\re A$
is the restriction of $f$ to $A$ and $f[A] :=
rng (f \re A)$ is the image of $A$ under $f$. $\forall^\infty n$ 
means {\it for all but finitely many $n$}, and $\exists^\infty n$
is used for {\it there are infinitely many $n$}. 

$\omoms$ ($\omloms$, respectively) denotes the infinite (finite, resp.)
subsets of $\omega$; $\omup$ ($\omlup$, resp.) stands for the
strictly increasing functions from $\omega$ to $\omega$ (for the
strictly increasing finite sequences of natural numbers, resp.).
Identifying subsets of $\omega$ with their increasing enumerations
naturally identifies $\omoms$ and $\omup$. We reserve letters like
$\sigma , \tau$ for elements of $\omlom$ and $\omlup$, and letters
like $s,t$ for elements of $\omloms$. $\ha$ is used for concatenation
of sequences (e.g., $\sigma\ha\la n\ra$). Given a tree $T\sub\omlom$, we
denote by $stem(T)$ its stem, and by $[T] := \{ f\in\omom; \; \forall n \; (f\re n
\in T )\}$ the set of its branches. Given $\sigma \in T$,
we let  $T_\sigma := \{ \tau\in T ; \;
\tau\sub
\sigma \;\lor\;\sigma\sub\tau\}$, the restriction of $T$
to $\sigma$, and $succ_T (\sigma) := \{ n \in\omega ; \;
\sigma \ha\la n\ra \in T \}$. For $A , B \sub
\omega$, we say $A\sub^* B$ ({\it $A$ is almost included in $B$})
iff $A\sem B$ is finite. If  $\A \sub \omoms$ and $B \in\omoms$
satisfies $B \sub^* A$ for all
$A\in\A$, we call $B$ a {\it pseudointersection} of $\A$.
A sequence $\T = \la T_\alpha ; \; \alpha < \kappa \ra$
is called a {\it $\kappa$--tower} (or: {\it tower of height
$\kappa$}) iff $T_\beta \sub^* T_\alpha$ for $\beta \geq \alpha$
and $\T$ has no pseudointersection.

Concerning forcing, let $\PP$ be a p.o. in the ground model $V$. $\PP$--names
are denoted by symbols like $\dot f$, $\dot X$, ..., and for their
interpretations in the generic extension $V[G]$, we use $f=
\dot f[G]$, $X=\dot X
[G]$... We often confuse Boolean--valued models $V^\PP$ and the
corresponding forcing extensions $V[G]$ where $G$ is
$\PP$--generic over $V$.
$\PP$ is called {\it $\sigma$--centered} iff there are $P_n \sub \PP$
with $\bigcup_n P_n = \PP$ and, for all $n$ and  $F\sub P_n$ finite,
there is $q\in\PP$ with $q\leq p$ for all $p\in F$.
$\star$ is used for two--step iteration (e.g. $\PP \star \dot \QQ$). 
If $\la \PP_\alpha , \dot \QQ_\alpha ; \; \alpha < \kappa \ra$ 
(where $\kappa$ is a limit ordinal) is an
iterated forcing construction with limit $\PP_\kappa$
(see [B] or [Je 1] for details) and  $G_\kappa$
is
$\PP_\kappa$--generic, we let $G_\alpha = G_\kappa \cap \PP_\alpha$ be the
restriction of the generic, and $V_\alpha = V[G_\alpha]
=V^{\PP_\alpha}$ stands for the
intermediate extension. In $V_\alpha$, $\PP_{[\alpha , \kappa)}$ denotes the
rest of the iteration.  $\CC_\kappa$ (where $\kappa$ is any ordinal)
stands for the p.o. adding $\kappa$ Cohen reals.
For sections 5 and 6, we assume familiarity with Easton forcing
(see [Je] or [Ku]) and the ways it can be factored. In particular,
we use that if $\PP$ is $ccc$ and $\QQ$ is $\omega_1$--closed 
(in $V$), then $\PP$ is still $ccc$ in $V^\QQ$ and $\QQ$ is
$\omega_1$--distributive in $V^\PP$. Recall that a p.o. $\QQ$
is called {\it $\lambda$--distributive} iff the intersection of fewer
than $\lambda$ open dense subsets of $\QQ$ is open dense.
In section 7, we shall need basic facts about club sets in $\omega_1$:
that the intersection of a countable family of clubs is
club, that given clubs $\{ C_\alpha ; \; \alpha < \omega_1 \}$, their
{\it diagonal intersection} $\{ \beta \in \omega_1 ; \; \forall \alpha < \beta
\; ( \beta \in C_\alpha ) \}$ is club, and that if $\PP$ is $ccc$ and $\forces_\PP ``\dot C$ is club",
then there is a club $D$ in the ground model such that $\forces_\PP
`` D \sub \dot C"$ (see [Ku, chapter II, $\S$ 6 and chapter VII, (H1)].

More notation will be introduced when needed.

\bigskip

{\bolds On the genesis of this paper and acknowledgements.} 
The first author is very much
indebted to the members of the logic group at Charles University, Praha: to
Bohuslav Balcar, Petr Simon and Egbert Th\"ummel for introducing him to the
world of characters and $\pi$--characters; to the latter for explaining him how
the cardinal characteristics of $r^0_\U$ could be read off from those of $\U$
in case $\U$ is a Ramsey ultrafilter. He gratefully acknowledges support
from the Center for Theoretical Study for his stay in
January/February 1995, and thanks Bohuslav
Balcar for having him invited. A preliminary version of this paper,
by the first author only, was circulated late in 1995. It consisted
of sections 2 to 4 and 7 of the present work and one more
section the results of which have been superseded.
Unfortunately, it contained several inaccuracies, and a few basic
results were not mentioned.

The main bulk of the important results in sections 5 and 6 were
proved by the second author in September 1996 while the first
author was visiting him at Rutgers University. Section 1 is joint
work. We thank Alan Dow, Martin Goldstern and Claude Laflamme for comments.
We also thank the referee for many valuable suggestions and for
detecting a gap in the original proof of Theorem 3.

\Bigskip


{\dunhg 1. Setting the stage --- some cardinal characteristics of ultrafilters}
\Smallskip

\no Let $\U$ be a non--principal ultrafilter on the natural numbers
$\omega$. We define the following four cardinal invariants
associated with $\U$.
\sm

\itemitem{$\pp (\U)$} $= \min \{ | \A | ; \; \A \sub \U \; \land \; \neg \exists
B\in\U \; \forall A \in\A \; (B \sub^* A ) \}$

\itemitem{$\pi\pp(\U)$} $ = \min \{ | \A | ; \; \A \sub \U \;\land\;
\neg \exists B \in\omoms \; \forall A \in \A \; (B \sub^* A) \}$

\itemitem{$\pi\chi (\U)$} $=\min \{ | \A | ; \; \A \sub \omoms \;\land\;\forall B
\in\U\;\exists A\in\A\; (A \sub^* B) \}$

\itemitem{$\chi(\U)$} $=\min \{ |\A |;\;\A\sub\U\;\land\;\forall B\in\U\;\exists
A\in\A\; (A\sub^* B)\}$
\sm

\no The definition of $\pp$ is dual to the one of $\chi$; similarly
$\pi\pp$ and $\pi\chi$ are dual. Therefore we can expect a strong
symmetry when studying these cardinals. Note that $\pp (\U) \geq
\omega_1$ is equivalent to saying {\it $\U$ is a $P$--point}. Ultrafilters
with $\pi\pp (\U) \geq\kappa$ are called {\it pseudo--$P_\kappa$--points
} in [Ny]. $\pi\chi(\U)$ is referred to as {\it $\pi$--character}, and
$\chi (\U)$ is known as the {\it character} of the ultrafilter
$\U$. Furthermore, a family $\A$ which has the property in the
definition of $\pi\chi(\U)$ ($\chi(\U)$, respectively) is called
a {\it $\pi$--base} ({\it base}, resp.) of $\U$. Both these
cardinals have been studied intensively, see e.g. [BK], [BS], [BlS],
[Ny] and [vM].


It is easy to see that for any ultrafilter $\U$, the following hold:
$\omega\leq\pp (\U) \leq \pi\pp(\U)$, $\pi\chi(\U)\leq
\chi(\U)\leq\cc$, and $\omega_1 \leq  \pi\pp (\U)$. Furthermore, $\pp (\U)$
is a regular cardinal, and we have $cf (\pi\chi(\U) ) \geq
\pp (\U)$. (The same holds with $\pi\chi$ replaced by $\chi$, see
Proposition 1.4 below for a stronger result.) To obtain more restrictions
on the possible values, and on the possible cofinalities,
of these cardinals, we need to introduce some
classical cardinal coefficients of the continuum. For $f,g\in\omom$, we
say {\it $g$ eventually dominates $f$}
($f \leq^* g$, in symbols) iff $f(n) \leq g(n)$ holds
for almost all $n\in\omega$. If $\U$ is an ultrafilter,
we say $g$ {\it $\U$--dominates} $f$ ($f\leq_\U g$,
in symbols) iff $\{ n ; \; f(n) \leq g(n) \} \in \U$.
\sm

\item{$\bb$} $= \min \{ | \F | ; \; \F \sub \omom \;\land\; \forall
g\in\omom\;\exists f\in\F\; (f \not\leq^* g ) \}$

\item{$\dd$} $=\min\{ |\F |; \;\F\sub\omom\;\land\;\forall g\in\omom\;
\exists f\in\F\; (g \leq^* f) \}$

\item{$\ss$} $=\min\{ | \A | ;\;\A \sub\omoms\;\land\;\forall B\in\omoms
\;\exists A\in\A \; (|A\cap B| = |(\omega\sem A ) \cap B| = \omega ) \}$

\item{$\rr$} $=\min\{ |\A | ; \; \A\sub\omoms\;\land\;\forall B\in\omoms
\;\exists A\in\A \; (A\sub^* B \;\lor\; A \sub^* \omega\sem B) \}$

\item{$\pp$} $=\min_\U \pi\pp(\U)$

\itemitem{\Cofu U} $=\min \{ |\F| ; \; \F\sub\omom\;\land\;
\forall g \in\omom \; \exists f \in \F \; (g \leq_\U f) \}$
\sm

\no $\bb$ and $\dd$ are dual, and so are $\ss$ and $\rr$. $\bb$ is called
{\it (un)bounding number}, $\dd$ is referred to as {\it dominating
number}, $\ss$ is known as {\it splitting number}, $\rr$ is called either
{\it reaping number} or {\it refinement number}, and $\pp$ is the
{\it pseudointersection number}. \Cofu {U} which
is self--dual is called the {\it cofinality}
of the ultraproduct $\omom / \U$. Families like $\F$ and $\A$
in the defining clauses of the first four of these numbers are referred to
as {\it unbounding}, {\it dominating}, {\it splitting}
and {\it reaping families}, respectively. It is known that $\pp$ and $\bb$ are
regular, that $\omega_1 \leq \pp \leq \bb \leq cf (\dd)$, that $\pp \leq \ss \leq \dd \leq
\cc$, and that $\bb \leq \rr\leq\cc$ (see [vD] and [Va]). Also recall that
$\pp = \cc$ is equivalent to $MA(\sigma$--centered) [Be], Martin's axiom
for $\sigma$--centered p.o.'s; thus all these cardinals equal $\cc$
under $MA$.

Concerning the relationship to the ultrafilter invariants, we see easily
that $\pi\pp (\U) \leq \ss$ and $\rr \leq \pi\chi(\U)$
for all ultrafilters $\U$. Also, $MA$ implies $\pi\pp (\U) =
 \cc$ for all $\U$, while there are 
(under $MA$) Ramsey ultrafilters
$\U$ with $\pp(\U) = \kappa$ for all regular $\omega_1 \leq \kappa \leq \cc$
[Lo, Th\'eor\`emes 3.9 et 3.12].
Furthermore, \Cofu {U} is regular and
$\bb \leq$ \Cofu U $\leq \dd$; for more results on
\Cofu {U} see [Bl], [Ca], [Ny], [SS] and the recent [BlM]. The following
proposition which relates the cofinality of $\omom / \U$ to other
invariants is well--known. We include a proof for completeness'
sake.
\sm

{\capit Proposition 1.1.} (Nyikos [Ny, Theorem 1 (i) and 3 (i)],
see also [Bl, Theorem 16]).

{\it (a) If $\pi\chi (\U) < \dd$, then \Cofu U $=\dd$.
Equivalently, $\max \{ \pi\chi (\U) ,$ \Cofu U$\} \geq\dd$.

(b) If $\pi\pp (\U) > \bb$, then \Cofu U $=\bb$.
Equivalently, $\min \{ \pi\pp (\U) ,$ \Cofu U$\} \leq \bb$.}
\sm

{\it Proof.} Given $f\in\omom$ and $A\in\omoms$ define
$f_A \in\omom$ by
$$f_A (n) : = \min \{ f(k) ; \; k \geq n \hbox{ and }k\in 
A\},$$
and note that if $g\in\omom$ is strictly increasing with
$g \leq_\U f$ then $g \leq^* f_A$ for any $A \sub^*
\{ n ; \; g(n) \leq f(n) \} \in \U$. $(\star)$

(a) If $\{ f^\alpha ; \; \alpha <$ \Cofu U$\}$ is cofinal
modulo $\U$ and $\{ A_\beta ; \; \beta < \pi\chi (\U) \}$
is a $\pi$--base, then $\{ f^\alpha_{A_\beta} ; \; \alpha
< $ \Cofu {U} and $\beta < \pi\chi(\U)\}$ is dominating by
$(\star)$.

(b) If $\kappa < \min \{ \pi\pp (\U) , $ \Cofu U$\}$
and $\{ g^\alpha ; \; \alpha < \kappa \} \sub \omom$
are strictly increasing, then find $f\in\omom$ with
$g^\alpha \leq_\U f$ for all $\alpha$. Put $A_\alpha
= \{ n ; \; g^\alpha (n) \leq f(n) \} \in \U$,
and find $A\sub^* A_\alpha$ for all $\alpha$. By $(\star)$,
we get $g^\alpha \leq^* f_A$ for all $\alpha$,
and the $g^\alpha$ are not unbounded. $\qed$
\sm

\no Since we always have $\pi\pp (\U) \leq \dd$ and $\pi\chi (\U)
\geq \bb$, we infer immediately
\sm

{\capit Corollary 1.2.} (Nyikos [Ny, Theorem 3 (viii)])
{\it For any ultrafilter $\U$, we have either $\pi\pp (\U)
\leq \bb$ or $\pi\chi (\U) \geq \dd$.} $\qed$
\sm

{\capit Corollary 1.3.} {\it $\pi\pp (\U) \leq \pi\chi (\U)$
holds for any ultrafilter $\U$.} $\qed$
\sm

\no We thus see that the four ultrafilter characteristics defined
at the beginning are, in fact, linearly ordered.

\bigskip

\ce{$\star\star\star$}

\bigskip

\no Unfortunately, we shall need some more ultrafilter coefficients
whose definition is not as nice as the one of the four above.
The reason for introducing these cardinals will become
clear in $\S\S$ 2 and 3. 
\sm

\itemitem{$\pp ' (\U) $} $=\min \{ |\A |; \; \A \sub \U \;\land\; \forall
\bar B \in [\U]^\omega \;\exists A\in\A\; \forall B\in\bar B \;
(B \not\sub^* A) \}$


\itemitem{$\pi\chi_\sigma (\U)$} $=\min \{ |\A | ;
\;\A\sub[\omega]^\omega\;\land\;
\forall \bar B \in [\U]^\omega \;\exists A\in\A\;\forall B\in\bar B \;
(A \sub^* B) \}$

\itemitem{$\chi_\sigma (\U)$} $=\min \{ |\A |; \;\A\sub [\U]^\omega \;\land\; \forall
\bar B \in [\U]^\omega \;\exists \bar A\in \A \;\forall B\in\bar B\;\exists
A\in\bar A \; (A \sub^* B) \}$
\sm

\no There is again some symmetry. For example, the cardinal which
is dual to $\pp ' (\U)$ can be defined as
\sm

\itemitem{$\chi ' (\U) $} $=\min \{ |\A| ; \; \A \sub [\U]^\omega \;\land\;
\forall B\in \U \;\exists \bar A \in \A \;\exists A\in\bar A \;
(A \sub^* B) \}$
\sm

\no Of course, we have $\chi ' (\U) = \chi (\U)$, and thus get nothing
new. Similarly, the primed version of $\pi\chi(\U)$,
as well as the $\sigma$--versions of $\pp(\U)$ and $\pi\pp (\U)$,
give us nothing new. One could define a primed version of $\pi\pp (\U)$,
but we won't need it.
Concerning the possible values of the primed
cardinal, we note that 
$\omega_1 \leq \pp ' (\U) \leq \pi\pp (\U) $ as well as $\pp (\U) \leq
\pp ' (\U) $. Furthermore,
$\pp' (\U)$ is regular, and we have the following result
which might be folklore:

\sm

{\capit Proposition 1.4.} {\it $cf(\chi (\U)) \geq \pp ' (\U)$.
In particular $\chi (\U)$ has uncountable cofinality.}
\sm

{\it Proof.} First note that if $\la \F_n ; \; n\in\omega\ra$
is a strictly increasing sequence of proper filters on $\omega$,
then $\F = \bigcup_n \F_n$ is not an ultrafilter. To see
this, choose a strictly decreasing sequence $\la A_n ; \; 
n\in\omega\ra$ of subsets of $\omega$ such that $A_0 = \omega$
and $A_{n+1} \in \F_{n+1} \sem \F_n$ for all $n$. Let
$B = \bigcup_n (A_{2n+1} \sem A_{2n+2})$ and 
$C = \bigcup_n (A_{2n} \sem A_{2n+1})$. Thus $B \cup C
=\omega$. Assume that $B \in \F$. Then $B \in \F_n$
for some $n$. Hence also $A_n \cap B \in \F_n$ and
$A_{n+1} \cap B \in \F_{n+1}$. If $n$ is even
we see $A_n \cap B \sub A_{n+1} \not\in \F_n$; if $n$ is
odd, we have $A_{n+1} \cap B \sub A_{n+2} \not\in \F_{n+1}$,
a contradiction in both cases. Therefore $B \notin \F$. Similarly we show
$C \notin \F$, and $\F$ is not an ultrafilter.

Now let $\kappa$ be regular uncountable and
assume $\la \F_\alpha ; \; \alpha < \kappa \ra$ is a strictly
increasing sequence of proper filters on $\omega$ with
$\F = \bigcup_\alpha \F_\alpha$.
Choose $A_{\alpha + 1} \in \F_{\alpha + 1} \sem \F_\alpha$.
Assume there are countably many $B_n \in\F$ such that
for all $\alpha$ there is $n$ with $B_n \sub A_{\alpha +1}$.
Then for some $\alpha_0 < \kappa$, $B_n \in \F_{\alpha_0}$ for all
$n$, a contradiction to the choice of $A_{\alpha_0+1}$.
Hence we see that $cf(\chi(\U)) \geq \pp '(\U)$ for any
ultrafilter $\U$. $\qed$

\sm

Also notice that $\pp (\U) = \pp ' (\U)$ iff $\U$ is
$P$--point. In particular, there are (in $ZFC$)
ultrafilters $\U$ with $\pp ' (\U)
> \pp (\U)$. Under $MA$ this can be strengthened to
\sm

{\capit Proposition 1.5.} (MA) {\it For each regular cardinal
$\kappa$ with $\omega_1 \leq \kappa \leq \cc$, there
is an ultrafilter $\U$ with $\pp (\U) = \omega$ and
$\pp' (\U) = \kappa$.}
\sm

{\it Proof.}  By Louveau's Theorem quoted above, there
is an ultrafilter $\V$ with $\pp (\V) = \kappa$.
Let $ X_n : = \{n\} \times\omega$ denote the vertical
strips. We define
an ultrafilter $\U$ on $\omega\times\omega$ by
\sm

\ce{$X \in \U \Loleriar \{ n ; \; \{ m ; \;
\la n,m\ra \in X \}  \in \V \} \in \V$.}
\sm

\no (We shall use again this type of construction in $\S$ 5.)
Note that the sets $Y_n : = \bigcup_{k\geq n} X_k$
witness $\pp (\U) = \omega$.

We are left with proving $\pp ' (\U) = \kappa$.
Given $A \in \U$, put $A_n = \{ m ; \; \la n,m \ra \in A \}$
and let $B_A = \{ n ; \; A_n \in
\V \} \in \V$. Notice that if $A \sub^* A'$ then
also $B_A \sub^* B_{A'}$. 

First take $\lambda < \kappa$
and let $\la A_\alpha ; \; \alpha < \lambda\ra$ be a sequence
from $\U$. By $\pp (\V) =\kappa$,
find $B \in \V$ with $B \sub^* B_{A_\alpha}$ for all
$\alpha$. Find $C_n \in \V$ such that $C_n \sub^* A_{\alpha,n}$
for all $\alpha$ with $A_{\alpha ,n} \in \V$.
Finally find $f\in\omom$ with $f (n) \geq \max (C_n \sem A_{\alpha ,
n})$ for almost all $n$ with $A_{\alpha ,n} \in \V$,
and all $\alpha$.
Then put $D_n = \bigcup_{k\geq n, k \in B} \{ k \}
\times (C_k \sem f(k))\in\U$.
It is now easy to check that for each  $\alpha < \lambda$
there is $n$ with $D_n \sub^* A_\alpha$. Hence $\pp' (\U)
\geq \kappa$.

Conversely, let $\la B_\alpha ; \; \alpha < \kappa \ra$ witness
$\pp  (\V) = \kappa$, and put $A_\alpha = \bigcup_{n\in B_\alpha}
X_n$. If we had $D_n \in \U$ such that for all $\alpha$ there
is $n$ with $D_n \sub^* A_\alpha$, then we would also get
$B_{D_n} \sub^* B_\alpha$, a contradiction. Thus
$\la A_\alpha ; \; \alpha < \kappa \ra$ witnesses
$\pp ' (\U) \leq \kappa$. $\qed$

\sm
\no On the other hand, it is easy to see that there is always an
ultrafilter $\U$ with $\pp ' (\U) = \omega_1$ (simply take
$\A = \{ A_\alpha ; \; \alpha < \omega_1 \}$ strictly $\sub^*$--decreasing,
let $\I$ be the ideal of pseudointersections of $\A$, and extend
$\A$ to an ultrafilter $\U$ with $\U \cap \I = \em$).
This should be seen as dual to the well--known fact (see e.g.
[vM, Theorem 4.4.2]) that
there is always an ultrafilter $\U$ with $\chi (\U) = \cc$.
\bigskip

\ce{$\star\star\star$}

\bigskip

\no To get more restrictions on the possible values of the
$\sigma$--versions of our ultrafilter characteristics, recall
the following cardinal invariants.

\sm

\item{$\rr_\sigma$} $=\min \{ | \A | ; \; \A \sub \omoms \;\land\;
\forall \bar B \in [\omoms]^\omega\;\exists A\in \A \;\forall B\in
\bar B \; (A\sub^* B\;\lor\; A\sub^*(\omega\sem B))\}$

\item{$\ppar$} $=\min \{ |\Pi| ; \; \Pi \sub 2^{[\omega]^2} \;\land\;
\forall A\in\omoms\;\exists\pi\in\Pi\;$ with $\pi[[A\sem n]^2] = 2$
for all $n\}$

\item{$\hhom$} $=\min \{ |\A| ;\; \A\sub\omoms\;\land\;$ for all
partitions $\pi : [\omega]^2 \to 2 \;$ there is $ A\in\A \;$ such that $A$ is
homogeneous for $\pi$ (that is, $|\pi[[A]^2]|=1)\}$

\sm

\no The partition cardinals $\ppar$ and $\hhom$ were introduced
by Blass [Bl 2, section 6]. It is known that $\ppar = \min \{ \ss , \bb \}$
and that $\hhom = \max \{ \rr_\sigma , \dd \}$ (see [Bl 2, Theorems 16 and
17], [Br, Proposition 4.2]). We see easily that  $\cc \geq
\chi_\sigma (\U) \geq \pi\chi_\sigma (\U)
\geq \rr_\sigma$, $\chi_\sigma (\U) \geq \chi(\U)$, $\pi
\chi_\sigma (\U) \geq \pi\chi (\U)$, $cf (\pi\chi_\sigma (\U)) \geq
\pp ' (\U)$, $cf (\chi_\sigma (\U)) \geq \pp ' (\U)$, and that
$\pi\chi_\sigma (\U) = \pi\chi(\U)$ as well as $\chi_\sigma (\U) =
\chi(\U)$ for $P$--points $\U$. We do not know whether $\chi_\sigma
(\U) > \chi (\U)$ is
consistent (see $\S$ 8 (1)), but we shall encounter
ultrafilters $\U$ with $\pi\chi_\sigma (\U) >
\pi\chi(\U)$ in section 5. 
The following proposition is simply a reformulation of the well--known
fact that Mathias forcing with a non--$P$--point adds a
dominating real. We include a proof for completeness' sake.
\sm

{\capit Proposition 1.6.} (Canjar, Nyikos, Ketonen, see [Ca 1, Lemma 4])
{\it Let $\U$ be an ultrafilter on $\omega$
which is not a $P$--point. Then: \par
(a) $\pi\pp (\U) \leq \bb$; \par
(b) $\pi\chi_\sigma (\U) \geq \dd$ and $\chi(\U) \geq \dd$.}
\sm

{\it Proof.} Let $\{ A_n ; \; n\in\omega\}\sub\U$ be decreasing 
with no infinite pseudointersection in $\U$; i.e. $A_{n+1}
\sub A_n$ and $| A_n \sem A_{n+1}| = \omega$ for all $n\in\omega$.
Given $f \in\omup$, let $A_f \in\U$ be such that $\min (A_f \cap
(A_n \sem A_{n+1})) \geq f(n)$ for all $n\in\omega$. Given $A
\in \omoms$, define $f_A (n) \in \omega $ by first finding the least 
$k \geq n$ with $A \cap (A_k \sem A_{k+1}) \neq \em$, if it exists,
and then putting $f_A (n) = \min (A \cap (A_k\sem A_{k+1}))$;
otherwise let $f_A (n) = 0$.

(a) Let $\kappa < \pi\pp (\U)$, $\{ f_\alpha ; \; \alpha < \kappa\}
\sub \omup$. Let $B$ be a pseudointersection of the family $\{
A_n ; \; n\in\omega \} \cup \{ A_{f_\alpha} ; \; \alpha < \kappa \}$.
It is easy to see that $f_B$ eventually dominates all $f_\alpha$.

(b) Let $\{ A_\alpha ; \; \alpha < \pi\chi_\sigma (\U) \}$ be a 
$\pi\sigma$--base
of $\U$. Given $f\in\omup$, let $\alpha$ be such that $A_\alpha \sub^* A_f
\cap A_n$ for all $n$.
Then $f_{A_\alpha}$ eventually dominates $f$. Thus $\{ f_{A_\alpha}
; \; \alpha < \pi\chi_\sigma (\U) \}$ is dominating. In case 
the $A_\alpha$ form a base, argue similarly: choose $\alpha$ such
that $A_\alpha \sub^* A_f$, etc.
$\qed$
\sm

\no We will see in 5.4 that $\pi\chi_\sigma$ and $\chi$
cannot be replaced by $\pi\chi$ in (b), in general.
We notice that the above result is also true for rapid ultrafilters
--- with an even easier argument. However, it may fail in general
(see the main results of [BlS] and [BlS 1]). The following proposition
has a flavor similar to Bartoszy\'nski's classical (and much more
intricate) result [Ba] that if \cov(measure) $\leq \bb$, then
\cov(measure) has uncountable cofinality. 
\sm

{\capit Proposition 1.7.} {\it If $\pi\pp (\U) \leq \bb$, then
$cf(\pi\pp (\U)) \geq\omega_1$.}
\sm

{\it Proof.} 
Assume $\lambda$ has countable cofinality and $\pi\pp (\U)
\geq\lambda$. We shall show $\pi\pp (\U) > \lambda $. 
Choose $\A \sub \U$ of size $\lambda$. Then $\A = \bigcup_n \A_n$
where $|\A_n| < \lambda$ and $\A_n \sub \A_{n+1}$. 
Hence we can find $X_n \in \omoms$
with $X_n \sub^* A$ for all $A\in \A_n$. For $A \in \A_n$
choose a function $f_A\in\omom$ with $X_k \sem A \sub f_A (k)$ 
for $k\geq n$. By assumption $\lambda < \bb$; hence there is
$f\in\omom$ with $f \geq^* f_A$ for all $A\in\A$.
Put $X := \{ \min (X_k \sem f (k)) ; \; k\in\omega\}$.
It's easy to check that $X \sub^* A$ for all $A\in\A$,
and we're done.
$\qed$
\sm

\no Proposition 1.6 and 1.7 together yield:

\sm

{\capit Corollary 1.8.} {\it If $\U$ is either not a $P$--point
or a rapid ultrafilter, then $\pi\pp (\U)$ has uncountable
cofinality.} $\qed$
\sm

For later use ($\S\S$ 2 and 3) we mention
the following characterization of $\chi_\sigma (\U)$.
\sm

{\capit Lemma 1.9.} {\it $\chi_\sigma (\U) = \min \{ | \A | ; \;
\A \sub \U^\omega \;\land\; \forall \la B_n ;\; n\in\omega\ra \sub\U
\;\exists \la A_n ; \; n\in\omega\ra\in\A\;\forall n\; (A_n \sub^* B_n) \}$.}
\sm

{\it Proof.} Denote the cardinal on the right--hand side by $\bar \chi_\sigma
(\U)$. $\chi_\sigma (\U) \leq \bar\chi_\sigma (\U)$ is trivial.
To see the converse, note that for $P$--points $\U$, both cardinals
coincide with the character. Hence assume $\U$ is not $P$--point;
then $\dd \leq \chi_\sigma (\U)$ by Proposition 1.6. Let $\{ f_\beta ; \;
\beta < \dd\}$ be a dominating family which is closed under
finite modifications (i.e. whenever $f\in\omom$ agrees with some
$f_\beta$ on all but finitely many places, then $f = f_\gamma$ for
some $\gamma < \dd$), and let $\{ \bar A_\alpha ; \;
\alpha < \chi_\sigma (\U) \}$ be a $\sigma$--base of $\U$. Let 
$\la A_{\alpha , n} ; \; n\in\omega \ra$ enumerate $\bar A_\alpha$;
without loss $A_{\alpha ,n+1} \sub^* A_{\alpha , n}$. Put $A_{\alpha ,
\beta , n} ' = A_{\alpha , f_\beta (n) }$; we leave it to the reader to verify
that $\{ \la A_{\alpha ,\beta,n} ' ; \; n\in\omega\ra ; \; \alpha
< \chi_\sigma (\U) , \beta < \dd \}$ satisfies the defining clause
of $\bar\chi_\sigma (\U)$. $\qed$

\vfill\eject


{\dunhg 2. Characterizations of the coefficients of the Ramsey ideal}
\Smallskip

\no Let $\I$ be a non--trivial ideal on the Baire space $\omom$ (or on one
of its homeomorphic copies, $\omoms$ or $\omup$) containing
all singletons. $\F\sub\I$ is a {\it base} of $\I$ iff given $A\in\I$
there is $B\in\F$ with $A\sub B$.
We introduce the following four cardinal invariants associated
with $\I$.
\sm

\itemitem{\Add I} $=\min\{ | \F |;\;\F\sub\I \;\land\; \bigcup \F
\not\in\I\}$

\itemitem{\Cov I} $=\min\{ | \F |;\;\F\sub\I\;\land\;\bigcup\F =\omom\}$

\itemitem{\Non I} $=\min\{ |F|;\; F\sub\omom\;\land\; F\not\in\I\}$

\itemitem{\Cof I} $=\min\{ |\F | ;\;\F\sub\I \;\land\; \F$ is a base of $\I\}$
\sm

\no These cardinals are referred to as {\it additivity}, {\it covering},
{\it uniformity} and {\it cofinality}, respectively. They
have been studied intensively in case $\I$ is either
the ideal of Lebesgue null sets or the ideal of meager
sets [BJ 1] and in some other cases as well. We note
that one always has \Add I $\leq$ \Cov I $\leq$ \Cof I{} and
\Add I $\leq$ \Non I $\leq$ \Cof I; furthermore, \Add I{} is regular,
and $cf($\Non I$)\geq $ \Add I, as well as $cf($\Cof I$)\geq $ \Add I.

Given an ultrafilter $\U$ on $\omega$, we define the {\it Mathias forcing}
associated with $\U$, $\MM_\U$ [Ma], as follows. Conditions are pairs
$(r,U)$ with $r\in \omloms$ and $U\in\U$ such that $\max (r) < \min (U)$.
We put $(s,V) \leq (r,U)$ iff $s \supseteq r$, $V \sub U$
and $s\sem r\sub U$. The Mathias p.o. is $\sigma$--centered and hence
$ccc$. It generically adds a real $m\in\omoms$ which is almost included
in all members of $\U$. For $(r,U) \in \MM_\U$, we
let $[r,U] = \{ A\in \omoms ; \; r \sub A \sub r \cup U \}$. 
The ideal of {\it nowhere Ramsey sets} with respect to $\U$
(or {\it Ramsey null sets}) 
$r^0_\U$ consists of all $X \sub \omoms$ such that given
$(r,U) \in \MM_\U$ there is $(r,V) \leq (r,U)$ with $X \cap
[r,V] = \em$. We notice that the connection between Mathias forcing
and the Ramsey ideal is like the one between Cohen (random, resp.)
forcing and the meager (null, resp.) ideal.

The main goal of this section is to characterize the four cardinal
coefficients introduced above for the ideal $r^0_\U$ in terms of
the cardinals in section 1. This extends a result of Louveau who had
already proved that the additivity of $r^0_\U$ coincides with
$\pp (\U)$.
For our characterizations we shall need  \sm

{\capit Lemma 2.1.} (Louveau, [Lo, Lemme 3.3]) {\it Let $\U$ be a 
$P$--point and $\phi : \omloms \to \U$. Then there is $U \in \U$ such
that $\{ s \in \omloms ; \; U \sem s \sub \phi (s) \}$ is cofinal
in $\omloms$.}
\sm

{\it Proof.} We include a proof to make the paper self--contained. Assume
$\U$ and $\phi$ are as required. Since $\U$ is a $P$--point, there is $U
\in\U$ with $U \sub^* \phi (s)$ for all $s\in \omloms$. Construct recursively
finite sets $A_i \sub U$ for $i\in\omega$ by putting $A_0 : = U
\sem \phi(\em)$ and $A_{i+1} := U \sem \bigcap \{ \phi (s) ; \;
\max(s) \leq \max(A_i) \}$. Then  either we have $U \sem \bigcup_i A_i
\in\U$, and this set is as required; or $\bigcup_i A_i \in \U$,
and one of the sets $\bigcup_i (A_{2i + 1} \sem A_{2i})$, $\bigcup_i
(A_{2i} \sem A_{2i - 1})$ lies in $\U$ and satisfies the conclusion
of the Lemma.
$\qed$
\bigskip


{\bolds Theorem 1.}  {\it Let $\U$ be an ultrafilter on $\omega$. Then:}
\par
(a) (Louveau [Lo, Th\'eor\`eme 3.7])
\add$(r^0_\U) =  \pp (\U)$; \par
(b) \cov$(r^0_\U) = \pi \pp (\U)$; \par
(c) \non$(r^0_\U) = \pi\chi (\U)$; \par
(d) \cof$(r^0_\U) = \chi_\sigma (\U)$.
\sm

In case $\U$ is a Ramsey ultrafilter, (a) through (d) were proved by Egbert
Th\"ummel. (Note that $\chi (\U) = \chi_\sigma (\U)$ in this case.)
\sm

{\it Proof.} Before plunging into the details, we describe natural
ways of assigning sets in the ideal to sets in the ultrafilter,
and vice--versa. Given $A \in \U$, let $X = X(A) := \{ B \in
\omoms ; \; B \not\sub^* A \}$ and note that $ X(A) = [\omega]^\omega
\sem \bigcup_{s
\in [\omega]^{<\omega}} [s,A\sem (\max(s)+1)] \in r^0_\U$. Conversely,
given $Y \in r^0_\U$, we can find a sequence $\la B_s \in \U ; \;
s\in [\omega]^{<\omega}\ra$ such that $B_s \sub \omega \sem (\max
(s)+1)$, $B_s \sub B_{t}$ for
$t \sub s$ and $Y \sub Y(\la B_s ; \; s\in \omloms \ra) : = 
[\omega]^\omega \sem \bigcup_s [s,B_s] \in r^0_\U$. Thus sets of the form
$Y(\la B_s \ra)$ form a base of the ideal $r^0_\U$,
and it suffices to deal with such sets in order to prove the Theorem.
We shall do this without further mention. Also, whenever dealing
with sequences $\la A_s \in \U ; \; s \in \omloms\ra$ we shall tacitly
assume that $A_s \sub \omega \sem (\max (s) + 1)$ and $A_s \sub A_t$ for
$t \sub s$.
We group dual results together.
\sm

{\sanse (a) and (d); the inequalities \add$(r^0_\U) \leq \pp(\U)$ and
\cof$(r^0_\U) \geq \chi_\sigma (\U)$.} 
Let $\{ A_\alpha ; \; \alpha <
\pp (\U) \} \sub \U$ be a witness for $\pp (\U)$. Let $X_\alpha =
X(A_\alpha)$. To see that
$\bigcup_\alpha X_\alpha \not\in r^0_\U$, fix $Y = Y(\la B_s 
\ra) \in r^0_\U$. There is $\alpha <\pp (\U)$ with
$y: = B_\em \sem A_\alpha$ being infinite. This means $y \in X_\alpha
\sem Y$, and we're done.

For the second inequality, notice that given
$\la A_t \in \U ; \; t\in\omloms\ra$ and $\la B_t \in\U ; \; t\in\omloms\ra$
with $A_\em$ coinfinite and $B_t \sem A_s$ infinite for some $s$
and all $t$,
we can construct $y \in Y(\la A_t \ra) \sem Y(\la B_t \ra)$ as follows:
choose $k > \max(s)$ such that $k\not\in A_\em$, and let
$y := s \cup \{ k \} \cup B_{s\cup\{k\}} \in [s\cup\{k\} , B_{s\cup\{k\}}]$;
then $y \not\in [t, A_t]$ for $t\sub s$ by the choice of $k$,
and $y \not\in [t,A_t]$ for $t\supseteq s$ by the properties of $A_s$.
\hskip 1truecm $(\star)$

Now let $\{ Y_\alpha ; \; \alpha <$ \cof$(r^0_\U)
\}$ be a base of the ideal $r^0_\U$. Without loss $Y_\alpha =
Y(\la B_{\alpha ,s } \ra)$ with all $B_{\alpha , s}
\in \U$. 
Fix $\bar A \in [\U]^\omega$; making its sets smaller, if necessary, we may
assume that $\bar A = \{ A_s ; \; s\in\omloms \}$ with $A_s \sub A_{t}$
for $t \sub s$ and $A_\em$ being coinfinite. 
Let $Y := Y(\la A_s \ra)$, and
choose $\alpha <$ \cof$(r^0_\U)$ with $Y \sub Y_\alpha$.
By $(\star)$ we get that for all $A \in \bar A$, there is
$s$ with $B_{\alpha , s} \sub^* A$, and we're done.
\sm

{\sanse The inequalities \add$(r^0_\U) \geq \pp (\U)$ and \cof$(r^0_\U)
\leq \chi_\sigma (\U)$.}  We distinguish
two cases.
First assume $\U$ is not a $P$--point. Then the first inequality is
trivial by $\pp (\U) = \omega$. Concerning the second, 
let $\{  \{ A_{\alpha,s} ; \;
s\in \omloms \}  ; \;
\alpha < \chi_\sigma (\U) \}$ be a $\sigma$--base of $\U$,
recall from
Proposition 1.6 that $\chi_\sigma (\U) \geq \dd$, let $\{ f_\beta : \omloms
\to \omega ; \; \beta < \dd \}$ be a dominating family which is
closed under finite modifications,
and put $Y_{\alpha , \beta } := Y(\la A_{\alpha,s} \sem f_\beta (s) 
\ra)$. We claim that $\{ Y_{\alpha , \beta} ; \;
\alpha < \chi_\sigma(\U) , \beta < \dd \}$ is a base of $r^0_\U$.
For, given $Y = Y(\la B_s \ra) \in r^0_\U$ with
$B_s \in \U$ for all $s$, we can find first 
(by Lemma 1.9) an $\alpha$ with
$A_{\alpha,s} \sub^* B_s$ for all $s$ and then a $\beta$
with $A_{\alpha,s} \sem f_\beta (s) \sub B_s$ for all $s$. This easily
entails $Y \sub Y_{\alpha,\beta}$.

Now suppose $\U$ is a $P$--point. Given $\la B_s \in \U ; \;
s\in \omloms\ra$ satisfying additionally $B_s = B_t$ for $s$
and $t$ with $\max (s) = \max (t)$ (and thus $B_s \sub B_t$ for
$s,t$ with $\max (t) \leq \max (s)$), as well as $A \in \U$ such
that $\{ s\in \omloms ; \; A\sem s \sub B_s \}$ is cofinal in
$\omloms$, we have $Y(\la B_s \ra) \sub X(A)$. 
$(\star\star)$ \hskip 0.4truecm To see this, fix
$s\in\omloms$, and take an arbitrary $y \in [s, A \sem (\max(s)+1) ]$.
Find $t \supseteq s$ with $A \sem t \sub B_t$. Letting $k := \max
(t) +1$, we get $y\sem k \sub A \sem k \sub B_t \sub B_{y\cap k}$
which entails $y \in [y\cap k , B_{y\cap k}]$.

Given $\kappa < \pp(\U)$ and $\{ Y_\alpha ; \; \alpha < \kappa\}
\sub r^0_\U$ where $Y_\alpha = Y(\la B_{\alpha,s} \ra)$ with
all $B_{\alpha,s} \in \U$, we find by Lemma 2.1 $A \in \U$
such that $\{ s\in\omloms ; \; A \sem s \sub B_{\alpha , s} \}$
is cofinal in $\omloms$ for all $\alpha$. Thus $\bigcup_\alpha Y_\alpha
\sub X(A) \in r^0_\U$ by $(\star\star)$. Dually,
if $\{ A_\alpha ; \; \alpha < \chi(\U) \}$ is a base of $\U$,
we claim that the sets $X_\alpha = X( A_{\alpha}
)$ form a base of our ideal. To see this,
take $Y = Y(\la B_s \ra) \in r^0_\U$ where
$B_s \in \U$.  By Lemma 2.1 find $\alpha
< \chi(\U)$ such that $\{ s\in\omloms ; \; A_\alpha\sem s 
\sub B_s \}$ is cofinal in $\omloms$, and conclude by $(\star\star)$.
\sm

{\sanse (b) and (c); the inequalities \cov$(r^0_\U) \leq \pi \pp (\U)$
and \non$(r^0_\U) \geq \pi\chi (\U)$.} This is easy. Given 
a witness $\{ A_\alpha \in \U ; \; \alpha < \pi\pp (\U) \}$
for $\pi \pp (\U)$, let $X_\alpha = X(A_\alpha)$. The $X_\alpha$
cover the reals, for, given $x\in \omoms$, there is $\alpha$ with $x \not\sub^*
A_\alpha$ which entails $x\in X_\alpha$. Dually, given
$\{ x_\alpha \in \omoms ; \; \alpha < $ \non$(r^0_\U) \} \not\in
r^0_\U$ and
$A \in \U$, there is $\alpha <$ \non$(r^0_\U)$ with $x_\alpha
\notin X(A)$ which means that $x_\alpha \sub^* A$. This shows that
the $x_\alpha$ form a $\pi$--base of $\U$.
\sm

{\sanse The inequalities \cov$(r^0_\U) \geq \pi \pp (\U)$ and
\non$(r^0_\U) \leq \pi\chi (\U)$.} We prove the second inequality first.
Let $\{ x_\alpha \in \omoms ; \; \alpha < \pi\chi (\U)\}$
be a $\pi$--base of $\U$. Given $n\in\omega$, let $x_{\alpha , n}
= x_\alpha \sem n$. We note that $\{ x_{\alpha , n} ; \; \alpha
< \pi\chi (\U) , n \in \omega \} \not\in r^0_\U$, because,
given $Y = Y(\la B_s \ra) \in r^0_\U$ with all $B_s \in \U$,
we find $\alpha$ with $x_\alpha \sub^* B_\em$ and thus $n\in\omega
$ with $x_{\alpha  , n} \sub B_\em$, that is $x_{\alpha , n} \not\in
Y$.

Next, let $\kappa < \pi\pp (\U)$ and  $\{ Y_\alpha ; \; \alpha
< \kappa\} \sub r^0_\U$; without loss $Y_\alpha
= Y (\la B_{\alpha , s} \ra)$ with $B_{\alpha,s}
\in\U$. We want to show that the $Y_\alpha$'s
do not cover the reals. We distinguish two cases.

First assume $\U$ is not a $P$--point. By assumption, we find
$x \in\omoms$ with $x \sub^* B_{\alpha, s}$ for all $\alpha$ and all
$s$. Define $g_\alpha : \omega \to x$ for $\alpha < \kappa$
recursively by:
$$\eqalign{g_\alpha (0) &:= \min \{ k ; \; x \sem k \sub B_{\alpha ,
\em} \} \cr g_\alpha (n+1) &:= \min
\{ k ; \; x\sem k \sub \bigcap_{s \sub g_\alpha (n) + 1} B_{\alpha , s} \}.
\cr } $$
By Proposition 1.6, we find $g : \omega \to x$ strictly increasing
and eventually dominating all
$g_\alpha$'s. Put $y : = rng (g)$. To complete the argument,
we shall show that $y \not\in \bigcup_\alpha Y_\alpha$. Fix $\alpha$.
Let $n_0$ be minimal with $g(n) \geq g_\alpha (n)$ for all $n \geq n_0$,
and put $s := \{ g(i) ; \; i < n_0 \}$. It is easy to see that
$y \in [s, B_{\alpha , s}]$.

Finally assume $\U$ is a $P$--point. By Lemma 2.1 we find $A_\alpha
\in\U$ such that $\{ s\in \omloms ; \; A_\alpha \sem s \sub
B_{\alpha ,s } \}$ is cofinal in $\omloms$. By assumption we find
$y \in \omoms$ with $y \sub^* A_\alpha$ for all $\alpha$. We show again
$y \not\in \bigcup_\alpha Y_\alpha$. Fix $\alpha$ and choose $s \in
\omloms$ with $y\sem s \sub A_\alpha$. Find $t \supseteq s$ with
$A_\alpha \sem t \sub B_{\alpha , t}$. Then, letting $k := \max(t) +1$,
we have $y \sem k \sub A_\alpha \sem k \sub B_{\alpha,t} \sub B_{\alpha,y\cap k}
$ which implies $y \in [y\cap k , B_{\alpha , y\cap k}]$.
This completes the proof of the Theorem.
$\qed$
\bigskip

{\capit Corollary 2.2.} {\it Let $\U$ be an ultrafilter
on $\omega$. Then:

(a) \add$(r^0_\U) \leq$ \cov$(r^0_\U) \leq$ \non$(r^0_\U) \leq
$ \cof$(r^0_\U)$. \par
(b) $\pp \leq $ \cov$(r^0_\U)$; in particular, $MA$ implies that
\cov$(r^0_\U) = \cc$. \par
(c) $cf($\cof$(r^0_\U)) \geq\omega_1$.}
\sm

{\it Proof.} All this follows from the Theorem and the results in $\S$
1, in particular Corollary 1.3 and Proposition 1.4. $\qed$
\sm

\no The fact that the cardinal coefficients of $ccc$--ideals of
the form $r^0_\U$ are linearly ordered distinguishes them from the 
$ccc$--ideals of meager and null sets (see [BJ 1]). Note,
however, that the cardinal coefficients of the closely related, but
non--$ccc$, ideal of {\it Ramsey null sets} [El] ({\it nowhere
Ramsey sets}) $r^0$ are also linearly
ordered. Namely, one has \add$(r^0 ) =$ \cov$(r^0) =\hh \leq $ \non$
(r^0) = \cc <$ \cof$(r^0)$ where $\hh$ is as usual the {\it
distributivity number} of $\P (\omega) / fin$ (this is due
to Plewik [Pl]).

Results with a flavor similar to our Theorem 1 were established
independently by Matet [M, section 10]. He considers a situation
which is both more general (filters on arbitrary regular $\kappa$
instead of ultrafilters on $\omega$) and more restricted (combinatorial
properties imposed on the filters) so that our results
are, to some extent, orthogonal.

For more results on the coefficients of $r^0_\U$,
see, in particular, Theorem 3, Theorem 4(c) and Corollary 6.5. 

\Bigskip


{\dunhg 3. Characterizations of the coefficients of the Louveau ideal}
\Smallskip

\no A tree $T \sub \omlup$ is called {\it Laver tree} iff for all
$\sigma \in T$ with
$\sigma \supseteq stem (T)$, the set $succ_T (\sigma) := \{ n \in \omega ; \; \sigma\ha\la
n\ra\in T\}$ is infinite. 
Given an ultrafilter $\U$ on $\omega$, we define the {\it Laver forcing}
associated with $\U$, $\LL_\U$ (see [Bl 1, section 5] or [JS, section 1]), as
follows. Conditions are Laver trees $T \sub\omlup$ such that for all $\sigma
\in T$ with $\sigma \supseteq stem (T)$, we have $succ_T (\sigma) \in \U$.
We put $S \leq T$ iff $S\sub T$; furthermore $S \leq_0 T$ iff
$S\leq T$ and $stem (S) = stem (T)$. The Laver p.o. is again a
$\sigma$--centered p.o. The {\it Louveau ideal} $\ell^0_\U$ consists
of all $X \sub \omup$ such that given $T \in \LL_\U$ there
is $S \leq T$ with $X \cap [S] = \em$ [Lo]. Louveau proved that
$\ell^0_\U$ is a $\sigma$--ideal and that there is a topology
$\G^\infty_\U$ on $\omup$, finer than the usual topology, such that
$\ell^0_\U$ is the ideal of the $\G^\infty_\U$ nowhere dense sets
which coincides with the $\G^\infty_\U$ meager sets [Lo, 1.11 et 1.12].
(This should be compared with Ellentuck's classical results [El] 
on nowhere Ramsey sets.)
He also showed that $\ell^0_\U = r^0_\U$ in case $\U$ is a Ramsey
ultrafilter [Lo, Propositions 1.3 et 3.1]. In the same vein, Blass [Bl 1, pp.
238--239] and Judah--Shelah  [JS, Theorem 1.20] observed that $\LL_\U$ and
$\MM_\U$ are forcing equivalent for Ramsey $\U$.

We are now heading for  characterizations of the cardinals
\add, \cov, \non{} and \cof{} for $\ell^0_\U$ in terms of the characteristics
of $\U$ introduced in section 1.
\bigskip

{\bolds Theorem 2.} {\it Let $\U$ be an ultrafilter on $\omega$. Then:}
\par
(a) \add$(\ell^0_\U) = \min \{ \pp ' (\U) , \bb \}$; \par
(b) \cov$(\ell^0_\U) = \min \{ \pi \pp (\U) , \bb \}$; \par
(c) \non$(\ell^0_\U) = \max \{ \pi\chi_\sigma (\U) , \dd \}$; \par
(d) \cof$(\ell^0_\U)= \max \{ \chi_\sigma (\U) , \dd \}$.
\sm

{\it Proof.} As in the proof of Theorem 1,
we shall stress the symmetry of the arguments,
and start by fixing some notation concerning the correspondence
between sets in $\U$ and sets in $\ell^0_\U$.
For $A \in \U$, let $X = X(A) : = \{ y \in \omup ; \; \exists^\infty
n \; (y(n) \notin A ) \}$; given additionally
$\sigma \in \omlup$, let $T=T_\sigma (A)$
be the Laver tree with stem $\sigma$ and $succ_T (\tau) = A \sem
(\tau (|\tau| - 1) + 1)$ for $\sigma\sub\tau\in T$. Then we have
$X(A) = \omup \sem \bigcup_{\sigma\in\omlup} [T_\sigma (A)]
\in \ell^0_\U$. Conversely, given $Y \in \ell^0_\U$,
we can find a sequence $\la B_\sigma \in \U ; \; \sigma \in \omlup \ra$
satisfying $B_\sigma \sub \omega \sem (\sigma (|\sigma| - 1) +1)$
and $B_\sigma \sub B_\tau$ for $\tau \sub \sigma$, and
such that $Y \sub Y(\la B_\sigma \ra )
:= \omup \sem \bigcup_\sigma [T_\sigma (\la B_\tau 
\ra )]$, where $T_\sigma (\la B_\tau 
\ra)$ is the Laver tree $T$ with stem $\sigma$
and $succ_T (\tau) = B_\tau$ for $\sigma\sub\tau\in T$. Again,
we use this convention without further comment.
\sm

{\sanse (a) and (d); the inequalities \add$(\ell^0_\U) \leq \pp ' (\U)$ and
\cof$(\ell^0_\U) \geq \chi_\sigma (\U)$.}
Notice first that given $\la A_\tau \in \U ; \; \tau \in
\omlup\ra$ and $\la B_\tau \in \U ; \; \tau \in \omlup \ra$,
if $y \in \omup$ satisfies $y(i) \in B_{y \re i} \sem
A_{y \re i}$ for almost all $i$, then $y \in Y(\la A_\tau \ra)
\sem Y(\la B_\tau \ra)$. \hskip 1truecm $(\star)$

Let $\{ A_\alpha ; \; \alpha < \pp ' (\U) \}$ be a witness for
$\pp ' (\U)$. We show that $\bigcup_\alpha X_\alpha \not\in \ell^0_\U$
where $X_\alpha = X(A_\alpha)$. This is easy, for, given
$Y=Y(\la B_\tau \ra)$ with $B_\tau \in \U$, we find $\alpha$
with $B_\tau \not\sub^* A_\alpha$ for all $\tau$, and then
construct $y\in \omup$ with $y(i) \in B_{y \re i} \sem A_\alpha$
for all $i$. This gives $X_\alpha \not\sub Y$ by $(\star)$.

Dually, let $\{ Y_\alpha ; \; \alpha < $ \cof$(\ell^0_\U) \}$
be a base of $\ell^0_\U$; without loss $Y_\alpha = Y(\la B_{\alpha ,
\sigma} \ra)$ with $B_{\alpha , \sigma} \in \U$. We claim that
the $\bar B_\alpha = \{ B_{\alpha , \sigma} ; \; \sigma \in
\omlup \}$ form a witness for $\chi_\sigma (\U)$. Let $\bar A
= \{ A_\sigma ; \; \sigma \in \omlup \}\in [\U]^\omega$ be given, and
find $\alpha$ with $Y:=Y(\la A_\sigma \ra) \sub Y_\alpha$. Assume
there were $\tau \in \omlup$ with $B_{\alpha , \sigma} \not\sub^*
A_\tau$ for all $\sigma$; then we could construct $y \in
\omup$ with $\tau \sub y$ and $y(i) \in B_{\alpha , y\re i} \sem A_{y \re i}$
for all $i\geq |\tau|$; this would contradict $Y\sub Y_\alpha$
by $(\star)$. Thus for all $A \in \bar A$ we find $B \in
\bar B_\alpha$ with $B \sub^* A$, and we're done.
\sm

{\sanse The inequalities \add$(\ell^0_\U) \geq \min \{ \pp ' (\U) ,
\bb \}$ and \cof$(\ell^0_\U) \leq \max \{ \chi_\sigma (\U) ,
\dd \}$.} First note that given $\la A_\tau \in\U ; \; \tau
\in \omlup \ra$ and $\la B_\tau \in \U ; \; \tau \in
\omlup \ra$ with $B_\tau \sub A_\tau$ for almost all
$\tau$, we have $Y(\la A_\tau \ra) \sub Y(\la B_\tau \ra)$.
$(\star\star)$ \hskip 0.4truecm
Let $\la \tau_n ; \; n\in\omega\ra$ enumerate $\omlup$.

Let $\kappa < \min \{ \pp' (\U) , \bb \}$, and let
$Y_\alpha  \in \ell^0_\U$ for $\alpha < \kappa$. Assume
$Y_\alpha = Y (\la A_{\alpha , \sigma } \ra)$ where $A_{\alpha
, \sigma} \in \U$. By $\kappa < \pp' (\U)$ find $\la B_n \in\U
; \; n\in\omega\ra$ such that for all $\alpha, \tau$,
there is $n$ with $B_n \sub^* A_{\alpha , \tau}$. Without loss,
$B_{n+1} \sub B_n$ for all $n$. Define $g_\alpha : \omega
\to \omega$ by:
$$g_\alpha (n) := \min \{ m ; \; B_m \sub^* A_{\alpha , \tau_n} \}.$$
Since $\kappa < \bb$, we can find $g \in \omom$ eventually
dominating all $g_\alpha$. Thus we have $B_{g(n)} \sub^* A_{\alpha,
\tau_n}$ for all $\alpha$ and almost all $n$. Define a function $h_\alpha $ 
for all $n$ with $B_{g(n)} \sub^* A_{\alpha , \tau_n}$ by:
$$h_\alpha (n) := \min \{ m ; \; B_{g(n)} \sem m \sub A_{\alpha ,
\tau_n} \}.$$
Let $h$ eventually dominate all $h_\alpha$. Then $B_{g(n)} \sem
h(n) \sub A_{\alpha , \tau_n }$ for all $\alpha$ and almost
all $n$. Put $B_{\tau_n} = B_{g(n)} \sem h(n)$. By $(\star\star)$
we have $Y_\alpha \sub Y(\la B_\tau \ra)$ for all $\alpha$,
and $\bigcup_\alpha Y_\alpha \in \ell^0_\U$ follows.

Dually, let $\{ \{ B_{\alpha , \tau} ; \; \tau\in
\omlup\} ;\; \alpha < \chi_\sigma
(\U)  \}$ be a $\sigma$--base of $\U$, and let
$\{ f_\alpha: \omlup\to\omega ; \; \alpha < \dd \}$ be a dominating
family. For $\alpha < \chi_\sigma (\U)$ and $\beta <
\dd$, let $Y_{\alpha , \beta}: = Y(\la B_{\alpha , \tau}
\sem f_\beta (\tau) \ra)$. Given $Y = Y(\la A_\tau \ra)
\in \ell^0_\U$ with $A_\tau \in \U$, first find 
(by Lemma 1.9) $\alpha$
with $B_{\alpha , \tau} \sub^* A_\tau$ for  all $\tau$,
and then $\beta$ with $B_{\alpha , \tau} \sem f_\beta (\tau)
\sub A_\tau$ for almost all $\tau$. By $(\star\star)$ we have
$Y \sub Y_{\alpha,\beta}$, and thus the $Y_{\alpha,\beta}$ form
a base of the ideal $\ell^0_\U$.
\sm

{\sanse (b) and (c); the inequalities \cov$(\ell^0_\U) \leq
\bb$ and \non$(\ell^0_\U) \geq \dd$.} Given $f\in\omom$ and
$\tau\in\omlup$ let $B_{f,\tau} : = \omega \sem \max\{ \tau(|\tau|
-1) + 1 , f(|\tau|) \}$ and $Y_f := Y(\la B_{f,\tau}\ra)$.
We easily see that, given an unbounded family $\{ f_\alpha\in\omom ; \;
\alpha < \bb \}$, we have $\bigcup_\alpha Y_{f_\alpha} = \omup$.
Dually, if $\{ f_\alpha \in\omup ; \;\alpha < $ \non$(\ell^0_\U) \}
\not\in\ell^0_\U$, then for each $f\in\omom$ there is $\alpha$
with $f_\alpha \not\in Y_f$ which means that $f \leq^* f_\alpha$;
hence the $f_\alpha$ form a dominating family. (Notice that this is just
a reformulation of the well--known fact that any Laver--like forcing
adds a dominating real.)
\sm

{\sanse The inequalities \cov$(\ell^0_\U) \leq \pi\pp (\U)$ and
\non$(\ell^0_\U) \geq \pi\chi_\sigma (\U)$.} Notice first that, given
a sequence $\la B_n \in \U ; \; n \in \omega \ra$ with $B_{n+1}
\sub B_n$ for all $n$, if we put $B_\sigma : = B_{|\sigma|}
\sem (\sigma(|\sigma| - 1) +1)$ and $Y: = Y(\la B_\sigma \ra)$,
then $y \in \omup \sem Y$ entails $rng (y) \sub^* B_n$ for all $n$.
\hskip 1truecm $(\dagger)$

Thus, given a witness $\{ A_\alpha ; \; \alpha < \pi\pp (\U) \}$
for $\pi\pp(\U)$, we must have $\bigcup_\alpha X ( A_\alpha ) = \omup$
(for, if $y$ were not in the union, $rng(y)$ would be a pseudointersection).
Similarly, if $\{ y_\alpha \in \omup ; \; \alpha < $ \non$(\ell^0_\U)
\} \not\in\ell^0_\U$, then for each $\la B_n \ra$ as in $(\dagger)$,
there is $\alpha$ with $y_\alpha \not\in Y(\la B_\sigma\ra)$, and
thus $y_\alpha \sub^* B_n$ for all $n$. This shows the $y_\alpha$ form
a $\pi\sigma$--base of $\U$.
\sm

{\sanse The inequalities \cov$(\ell^0_\U) \geq \min \{ \pi\pp
(\U) , \bb \}$ and \non$(\ell^0_\U) \leq \max \{ \pi\chi_\sigma
(\U) , \dd \}$.}
This is quite similar to the other two inequalities involving
$\min$ and $\max$ (see above). Given $\la A_\tau \in\U ; \; \tau\in
\omlup\ra$ and $\la B_\tau \in \omoms ; \; \tau\in\omlup\ra$ with $B_\tau
\sub A_\tau$ for almost all $\tau$, we have that any real $y\in\omup$
with $y(i) \in B_{y\re i}$ for all $i$ does not lie in $Y(\la A_\tau\ra)$.
\hskip 1truecm $(\ddagger)$

Let $\kappa < \min \{ \pi\pp (\U) , \bb \}$, and let $Y_\alpha
= Y(\la A_{\alpha,\sigma}\ra) \in \ell^0_\U$ with $A_{\alpha , \sigma }
\in \U$ for $\alpha < \kappa$. First find $B \in\omoms$
with $B \sub^* A_{\alpha,\sigma}$ for all $\alpha $ and $\sigma$,
then define $g_\alpha : \omlup\to\omega$ by:
$$g_\alpha (\tau) := \min \{ m ; \; B\sem m \sub A_{\alpha , \tau} \}.$$
Let $g: \omlup\to\omega$ eventually dominate all $g_\alpha$,
and put $B_\tau : = B \sem g(\tau)$. By $(\ddagger)$ we can construct
a real not in $\bigcup_\alpha Y_\alpha$, and the family we started with
is not a covering family.

Dually, let $\{ B_\alpha ; \; \alpha < \pi\chi_\sigma (\U) \}$ be $\pi\sigma$--base
of $\U$, and let $\{ f_\alpha : \omlup \to \omega ; \; \alpha < \dd \}$
be a dominating family. For $\alpha < \pi\chi_\sigma (\U)$ and $\beta
< \dd$ choose a real $y = y_{\alpha,\beta} \in \omup$ with $y(i) \in B_\alpha \sem
f_\beta (y\re i)$. Given $Y = Y (\la A_\tau\ra) \in \ell^0_\U$,
first find $\alpha$ with $B_\alpha \sub^* A_\tau$ for all $\tau$,
and then $\beta$ with $B_\alpha \sem f_\beta (\tau) \sub A_\tau$ 
for almost all $\tau$. By $(\ddagger)$ we know that $y_{\alpha,\beta}
\not\in Y$, and thus $\{ y_{\alpha,\beta} ; \; \alpha < \pi\chi_\sigma
(\U) , \beta < \dd \} \not\in \ell^0_\U$. This concludes the
proof of the Theorem.
$\qed$
\bigskip

{\capit Corollary 3.1.} {\it Let $\U$ be an ultrafilter on $\omega$. Then:
\par
(a) \add$(\ell^0_\U) \leq$ \cov$(\ell^0_\U) \leq$ \non$(\ell^0_\U) \leq
$ \cof$(\ell^0_\U)$. \par
(b) $\pp \leq $ \cov$(\ell^0_\U)$; in particular, $MA$ implies that
\cov$(\ell^0_\U) = \cc$. \par
(c) \cov$(\ell^0_\U) \leq \ppar$ and \non$(\ell^0_\U) \geq\hhom$.

(d) $cf($\cov$(\ell^0_\U)) \geq\omega_1$.}
\sm

{\it Proof.} All this is immediate from the Theorem and the
results concerning $\bb , \dd , \pp , \ss , \rr_\sigma , \ppar$
and $\hhom$ mentioned in $\S$ 1, in particular Corollary 1.3
and Proposition 1.7. $\qed$
\sm

\no Note that, since $\ell^0_\U$ is a $\sigma$--ideal,
both \non$(\ell^0_\U)$ and \cof$(\ell^0_\U)$ 
necessarily have uncountable cofinality.
Again, the cardinal coefficients are linearly ordered, like those
for the related non--$ccc$ {\it Laver ideal} $\ell^0$ (see
[GRSS]) --- or those for the $ccc$--ideal of meager sets in the
{\it dominating topology} (see [LR]); the latter topology in fact
sits strictly in between the standard topology on $\omup$ and
Louveau's topology $\G^\infty_\U$ which is relevant here.

Distinguishing the two cases whether or not $\U$ is a $P$--point,
we get somewhat nicer characterizations, by 1.6 and other
remarks in $\S$ 1.
\sm

{\capit Corollary 3.2.} {\it Assume $\U$ is a $P$--point. Then:}
\par
(a) \add$(\ell^0_\U) = \min \{ \pp (\U) , \bb \}$;
\par
(b) \cov$(\ell^0_\U) = \min \{ \pi\pp (\U) , \bb \}$;
\par
(c) \non$(\ell^0_\U) = \max \{ \pi\chi (\U) , \dd \}$;
\par
(d) \cof$(\ell^0_\U) = \max \{ \chi(\U) , \dd \}$.
\sm
\no{\it Assume $\U$ is not a $P$--point. Then:}
\par
(a) \add$(\ell^0_\U) = \pp ' (\U)$;
\par
(b) \cov$(\ell^0_\U) = \pi\pp (\U)$;
\par
(c) \non$(\ell^0_\U) = \pi\chi_\sigma (\U)$;
\par
(d) \cof$(\ell^0_\U) = \chi_\sigma (\U)$. $\qed$
\sm

\no For rapid $P$--points, the formulae get still simpler, and, in fact,
the invariants for the Ramsey ideal and the Louveau ideal coincide.
In view of Louveau's $r^0_\U = \ell^0_\U$ for Ramsey ultrafilters $\U$,
Theorem 2 provides an alternative way for calculating the
coeffients of $r^0_\U$.


We close this section with a diagram showing the relations between
the cardinal invariants considered in this work.
\bigskip

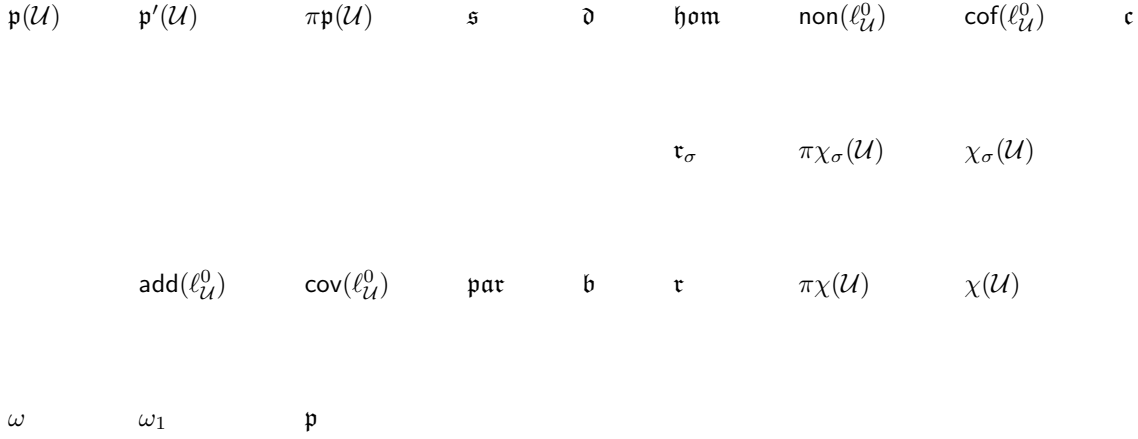
\halign{#  &\Hoskip # &\Hoskip  # &\Hoskip  # &\Hoskip  # &\Hoskip  # 
&\Hoskip # &\Hoskip # &\Hoskip # \cr
$\pp(\U)$ & $\pp' (\U)$ & $\pi\pp (\U)$ & $\ss$ & $\dd$ &
$\hhom$ & \non$(\ell_\U^0)$ & \cof$(\ell^0_\U)$ & $\cc$ \cr
\noalign{\Smallskip}
&&&&& $\rr_\sigma$ & $\pi\chi_\sigma(\U)$ & $\chi_\sigma (\U)$ & \cr
\noalign{\Smallskip}
& \add$(\ell^0_\U)$ & \cov$(\ell^0_\U)$ & $\ppar$ & $\bb$ & $\rr$ &
$\pi\chi (\U)$ & $\chi(\U)$ & \cr
\noalign{\Smallskip}
$\omega$ & $\omega_1$ & $\pp$ &&&&&& \cr}

\bigskip
\no Cardinals get larger when one moves up and to the right.
Dotted lines around three cardinals say that one of them is the
minimum or the maximum of the others for any ultrafilter $\U$
(which one of these alternatives holds being clear from the context).
For ease of reading, we omitted the inequality
$\pi\pp(\U) \leq
\pi\chi (\U)$.

\Bigskip


{\dunhg 4. Distinguishing the coefficients}
\Smallskip

\no Let $\U$ be a Ramsey ultrafilter. Since, by previous results,
\sm

\ce{\cov$(r^0_\U) = \pi\pp (\U) \leq \bb \leq \dd \leq \pi\chi (\U)
=$ \non$(r^0_\U)$,}
\sm

\no we can easily get the consistency of {\it there is a Ramsey ultrafilter
$\U$ with \cov$(r^0_\U) < $ \non$(r^0_\U)$.} For example, this can
be achieved by adjoining $\omega_2$ Cohen reals to a model of $CH$.
Furthermore, as mentioned in $\S$ 1, Louveau 
proved that $MA$ entails the
existence of a Ramsey ultrafilter $\U$ with
\sm

\ce{\add$(r^0_\U) = \pp(\U) = \omega_1 < \cc = \pi\pp (\U) =$
\cov$(r^0_\U)$.}
\sm

\no We complete this cycle of results by showing the
remaining consistency:

\bigskip

{\bolds Theorem 3.} {\it It is consistent with $ZFC$ that there is
a Ramsey ultrafilter $\U$ with \non$(r^0_\U) = \pi\chi (\U) <
\chi (\U) =$ \cof$(r^0_\U)$.} 
\sm

This answers half of the question in [Br, subsection 4.1].
As suggested by the referee, we note that this consistency has 
been well--known if one doesn't insist on $\U$'s Ramseyness (to see this
either use the Goldstern--Shelah model [GS] showing the consistency of
$\rr < \uu := \min_\U \chi (\U)$ and appeal to Proposition 7.1
below, or use
the Bell--Kunen model [BK] (cf. Remark 4.2 below)).

\sm

{\it Proof.} We start with a model $V$ of $CH$ and perform a finite
support iteration $\la \PP_\alpha , \dot\QQ_\alpha ; \; \alpha < \omega_1
\ra$ of $ccc$ p.o.'s. 
We build up the Ramsey ultrafilter $\U $
of $V_{\omega_1}$ along the iteration as a tower
of ultrafilters; in stage $2\cdot \alpha
+ 1$, we shall have the Ramsey ultrafilter $\U_{2\cdot\alpha + 1}$
in the model $V_{2\cdot\alpha + 1}$.
The details are as follows.
\sm

{\sanse Stage $\alpha$, $\alpha$ odd.}  In $V_\alpha$, we let $\QQ_\alpha
= \MM_{\U_\alpha}$ (Mathias forcing with the Ramsey
ultrafilter $\U_\alpha$, see $\S$ 2 for details). Denote by
$m_\alpha \in \omoms \cap V_{\alpha +1}$ the generic Mathias real
(which satisfies $m_\alpha \sub^* U$ for all $U \in \U_\alpha$).
\sm

{\sanse Stage $\alpha$, $\alpha$ even.} In $V_\alpha$, we let $\QQ_\alpha
 = \CC_{\omega_2}$, the p.o. for adding $\omega_2$ Cohen
reals. In $V_{\alpha + 1}$, we use the $\omega_2$ Cohen reals
to produce  the Ramsey ultrafilter
$\U_{\alpha +1} $ which extends either the filter $\F_\alpha$ generated
by $\bigcup_{\gamma < \alpha} \U_{2 \cdot \gamma +1}$ (in case
$\alpha$ is a limit) or the filter $\F_\alpha$  generated by $\U_{\alpha - 1} $
and $\omega \sem m_{\alpha - 1}$ (in case $\alpha $ is successor)
or the cofinite filter $\F_0$ (in case $\alpha = 0$).
This is
a standard construction (see, e.g., [Ca 2, Theorem 2], [BJ, $\S$ 3] or
[St, Theorem 5.2]) which we repeat
to make later arguments more transparent.

Let $V_{\alpha , \beta}, \; \beta \leq \omega_2,$ denote the model obtained
by adding $\beta$ of the Cohen reals (so $V_{\alpha , 0} = V_\alpha$
and $V_{\alpha , \omega_2} = V_{\alpha + 1}$). In $V_{\alpha + 1}$
enumerate the partitions of $\omega$ into infinite
subsets as $\la \la X_{\alpha,\beta,n} ; \; n\in\omega \ra ; \;
\beta < \omega_2\ra$ such
that $\la X_{\alpha , \beta,n };\; n\in\omega\ra \in V_{\alpha , \beta}$. 
Let $\U_{\alpha , 0}$ be a {\it careful} extension of $\F_\alpha$
to an ultrafilter of $V_{\alpha , 0}$ ({\it careful} will be
defined  later). Fix $\beta \leq \omega_2$,
and assume $\U_{\alpha , \gamma}$, a tower of ultrafilters in the respective
models $V_{\alpha , \gamma}$,
has been constructed. In case $cf(\beta) > \omega$, we let
$\U_{\alpha , \beta} = \bigcup_{\gamma < \beta} \U_{\alpha , \gamma}
$; and in case $cf(\beta) = \omega$, $\U_{\alpha , \beta}$ is a {\it
careful} extension of $\bigcup_{\gamma  < \beta} \U_{\alpha , \gamma}$
to the model $V_{\alpha , \beta}$.
In case $\beta = \gamma + 1$ do the following.

If $\bigcup_{k<n} X_{\alpha,\gamma ,k}
 \in \U_{\alpha , \gamma}$ for some $n$,  
we think of Cohen forcing $\CC$ as adjoining a subset of $\omega$,
called $c_{\alpha , \gamma}$, in the usual way. Call such $c_{\alpha , \gamma}$
{\it of the first kind}.
Otherwise think of Cohen forcing $\CC$ as adjoining a subset
$c_{\alpha , \gamma}$ of $\omega$ with $| c_{\alpha , \gamma}
\cap X_{\alpha,\gamma, n} | = 1$ for all $n$. 
Call such $c_{\alpha, \gamma}$ {\it of the second kind}. Then, 
by genericity, $c_{\alpha,
\gamma} \cap U$ is infinite for all $U \in \U_{\alpha,\gamma}
$. In both cases, let $\U_{\alpha,\beta}$ be a {\it careful}
extension of the filter generated by $\U_{\alpha,\gamma} $ and
$c_{\alpha,\gamma}$ to an ultrafilter in $V_{\alpha,\beta}$.
\sm

This completes the construction of the Ramsey ultrafilter
$\U = \bigcup_{\alpha < \omega_1} \U_{2\cdot\alpha + 1}$
in the resulting model $V_{\omega_1}$. The
$\omega_1$ Mathias reals $m_\alpha$ witness $\pi\chi (\U)
=\omega_1$. Thus we are left with showing $\chi(\U) = \omega_2$.
One inequality is clear because $\cc = \omega_2$. To see the other
one, we have to make precise what we mean by {\it careful}.

Let $\I$ be the family of subsets $I$ of $\omega$ for which
we can find pairwise disjoint finite sets $F_\alpha \sub \omega_1
\times \omega_2$,  $\alpha < \omega_1$, such that
$$\forall \alpha < \omega_1 \hskip 1.5truecm I \sub^* \bigcup_{(\beta,
\gamma) \in F_\alpha} c_{\beta , \gamma},$$
where the $c_{\beta,\gamma}$ denote the Cohen reals as explained
above. Clearly, $\I$ is an ideal (in $V_{\omega_1}$). Let  $\I_{\alpha , \beta } = \I 
\cap V_{\alpha ,\beta}$
for even $\alpha < \omega_1 $ and $\beta \leq \omega_2$. Note that
the above definition of $\I$ also makes sense, with the obvious adjustments, 
in each model $V_{\alpha , \beta}$. Call the resulting ideal 
$\I^{V_{\alpha , \beta}}$. Then one has $\I^{V_{\alpha , \beta}} =
\I_{\alpha , \beta}$, and thus $\I_{\alpha , \beta} \in V_{\alpha , \beta}$.

(The inclusion $``\sub"$ is obvious. To see $``\supseteq"$, note
first that if $I \in V_{\alpha , \beta}$ and $I \sub^* \bigcup_{(\gamma,
\delta) \in F} c_{\gamma , \delta}$, then $I \sub^* \bigcup_{(\gamma
,\delta) \in F'} c_{\gamma , \delta}$ where $F ' = F \cap (( \alpha
\times \omega_2) \cup (\{ \alpha \} \times \beta))$, by Cohenness of the
$c_{\gamma , \delta}$. Now assume (still in $V_{\alpha , \beta}$) that
$p \forces_{[(\alpha , \beta) , \omega_1 )} `` I \in \dot \I$ as witnessed
by $\dot F_\zeta , \zeta < \omega_1"$.
Find $p_\zeta \leq p$ (in $V_{\alpha , \beta}$) such that $p_\zeta$
decides $\dot F_\zeta$, say $p_\zeta \forces_{[(\alpha , \beta) , \omega_1)
} \dot F_\zeta = F_\zeta$. By the previous remark we may assume $F_\zeta
\sub (\alpha \times \omega_2 ) \cup (\{ \alpha \} \times \beta )$. Since
all the factors of the iteration satisfy Knaster's condition,
so does the quotient $\PP_{[(\alpha , \beta) , \omega_1 )}$. Thus
we may assume without loss that the $p_\zeta$ are pairwise compatible. 
Hence the $F_\zeta$ are pairwise disjoint. Therefore $I \in \I^{V_{\alpha
, \beta}}$ as required.)

We shall guarantee while extending the ultrafilter that 
$$(\star) \;\;\;\;\;  \U_{\alpha,\beta} \cap \I_{\alpha,\beta} = \em
\;\;\; \hbox{ for all even } \alpha < \omega_1 \hbox{ and all } \beta \leq \omega_2.$$
Such extensions will be called {\it careful}.
We have various cases to consider to see that this can be done.
\sm

\item{(1)} {\sanse Successor step.} Assume that $\U_{\alpha , \beta}
\cap \I_{\alpha ,\beta} = \em$ (where $\alpha < \omega_1$ is even and
$\beta < \omega_2$).
Since $V_{\alpha,\beta +1}$ is an extension by one Cohen real, we certainly
have $\la \U_{\alpha ,\beta } \ra \cap \I_{\alpha , \beta + 1} =
\em$ where $\la \U_{\alpha , \beta} \ra$ denotes the filter generated by
$\U_{\alpha , \beta}$ in $V_{\alpha ,\beta +1}$.
Next notice that $c_{\alpha , \beta} \cap U \not\in \I_{\alpha , \beta +1}
$ for all $U \in \U_{\alpha , \beta}$.  $(\dagger)$ \hskip 0.3truecm
To see this, fix such $U$. Assume $c_{\alpha , \beta}$
is of the second kind. Put $U_n = U \sem \bigcup_{k<n}
X_{\alpha,\beta,k}$. All $U_n$ lie in $\U_{\alpha , \beta}$ and thus
not in $\I_{\alpha , \beta + 1}$. Now note that whenever $F \sub 
(\alpha \times \omega_2) \cup (\{ \alpha \} \times \beta)$
is finite with $c_{\alpha , \beta} \cap U \sub^* \bigcup_{(\gamma , \delta)
\in F} c_{\gamma , \delta}$, then by Cohenness $U_n \sub^*
\bigcup_{(\gamma , \delta) \in F} c_{\gamma , \delta}$ for some $n$.
If $c_{\alpha , \beta}$ is of the first kind, the argument is even easier.
This shows $(\dagger)$. Now we can easily extend the filter
generated by $\U_{\alpha , \beta}$ and $c_{\alpha , \beta}$ to
an ultrafilter $\U_{\alpha , \beta+1}$ such that $\U_{\alpha , \beta +
1} \cap \I_{\alpha , \beta +1} = \em$.

\item{(2)} {\sanse Limit step.} Let $\alpha$ be even, and $\beta$ a limit ordinal.
If $\beta > 0$, assume that $\U_{\alpha,\gamma} \cap \I_{\alpha,\gamma}
= \em$ for $\gamma < \beta$. In case $cf(\beta) > \omega$ we get
$\U_{\alpha , \beta} \cap \I_{\alpha , \beta} = \em$ 
because any $U \in \U_{\alpha,\beta}$
lies already in $V_{\alpha,\gamma}$ for some $\gamma < \beta$, and thus cannot
be almost contained in the union of finitely many Cohen reals added
at a later stage. 
If $\beta = 0$ and $\alpha$ non--limit, it remains to see that
$\omega\sem m_{\alpha - 1} \cap U \not\in \I_{\alpha , 0}$
for all $U \in \U_{\alpha - 1}$ which is similar to, but easier
than, the argument for the $c_{\gamma , \delta}$ in (1). Whether or not
$\alpha$ is limit and whether $cf(\beta) = \omega$ or $\beta = 0$,
we can extend the given filter easily to $\U_{\alpha,\beta}$ such
that $\U_{\alpha,\beta} \cap \I_{\alpha,\beta} = \em$.
\sm

\no (1) and (2) clearly entail $(\star)$. If we had $\chi(\U)
\leq\omega_1$, we could find $U \in \U$ which is almost included
in $\omega_2$ of the Cohen reals which we added to $\U$ in the course
of the construction, and this would contradict $(\star)$. Thus
$\chi(\U) \geq \omega_2$, and the proof is complete.
$\qed$
\bigskip

{\capit Remark 4.1.} $\omega_1$ and $\omega_2$ in the above proof can be
replaced by arbitrary regular $\kappa < \lambda$. The argument is the
same: $\chi(\U) \leq \lambda$ by $\cc = \lambda$, $\chi(\U) \geq\lambda$
by $(\star)$, $\pi\chi(\U) \leq \kappa$ by the $\kappa$ Mathias
reals and $\pi\chi(\U) \geq\kappa$ by the fact that the iteration
has length $\kappa$ which implies $\pi\chi(\U) \geq\dd\geq$ \cov(meager)
$\geq	\kappa$ (where the first two inequalitites are in $ZFC$).
Here, \cov(meager) denotes the smallest size of a family of meager
sets covering the reals. It is well--known (and easy to see) that
$\dd \geq$ \cov(meager).
\sm

{\capit Remark 4.2.} The method of the proof of Theorem 3 can also
be used to show there is a Ramsey ultrafilter $\U$ with $\chi (\U)
=\cc$ in the Bell--Kunen model (see [BK]; this model is gotten
by a finite support iteration of $ccc$
p.o.'s of length $\omega_{\omega_1}$ over
a model of $CH$, forcing $MA$ for small p.o.'s at limit steps of the
form $\omega_{\alpha + 1}$; it satisfies $\cc = \omega_{\omega_1}$ and
$\pi\chi (\V) = \omega_1$ for all ultrafilters $\V$).
Hence it is consistent there is a Ramsey ultrafilter $\U$
with $\pi\chi (\U) = \omega_1 < \omega_{\omega_1} = \chi (\U)$.

\Bigskip
\vfill\eject


{\dunhg 5. A plethora of $\pi$--characters}
\Smallskip

\no This section is devoted to understanding the spectrum
of possible values for the $\pi$--character and clearing up
the relationship between $\pi$--character and $\pi\sigma$--character.
For this, we need to
discuss two ultrafilter constructions. 
First,  let
$\V$ and $\V_n$ be
ultrafilters on $\omega$, and define an ultrafilter $\U$ 
on $\omega\times\omega$ by
\sm

\ce{$X \in \U \Loleriar \{ n ; \; \{m; \; \la n,m \ra 
\in X \} \in \V_n \} \in \V$.}
\sm

\no (Note that we used this construction already once
in the proof of Proposition 1.5.) 
We call $\U$ the {\it $\V$--sum} of the $\V_n$, $\U = \sum_\V \V_n$.
In case all $\V_n$ are the same ultrafilter $\W$, we write
$\U = \V \times \W$ and call it the {\it product} of $\V$ and
$\W$. 
Then we have:
\sm

{\capit Proposition 5.1.} {\it
(a) $\min \{ \pi\chi (\V) , \sum_\V \pi\chi (\V_n )  \}
\leq
\pi\chi (\U) \leq \sum_\V  \pi\chi (\V_n ) $.

(b) If $\U = \V\times\W$, we have $\pi\chi (\U) = \pi\chi (\W)$.

(c) $
\pi\chi_\sigma (\U) \geq \pi\chi_\sigma (\V)$.

(d) If $\U = \V \times\W$, we have $\pi\chi_\sigma (\U) = \max
\{ \pi\chi_\sigma (\V) , \pi\chi_\sigma (\W) , \dd \}$.}
\sm

Here, given cardinals $\lambda_\alpha$, $\alpha \in R$, and
an ultrafilter $\D$ on $R$, $\sum_\D \lambda_\alpha$
denotes the {\it $\D$--limit} of the $\lambda_\alpha$, that is
the least cardinal $\kappa$ such that $\{ \alpha ; \;
\lambda_\alpha \leq \kappa \} \in \D$. \sm

{\it Proof.} For the purposes of this proof,
let $X_n  = \{ n \} \times \omega$
denote the vertical strips.

(a) The second inequality is easy, for we can take  the union
of $\pi$--bases of the appropriate $\V_n$'s (considered as ultrafilters
on the $X_n$'s) as a $\pi$--base
for $\U$. For the first inequality, let $\kappa < \min$
and $\A = \{ A_\alpha ; \; \alpha < \kappa \} \sub [\omega\times
\omega]^\omega$. 
We want to show $\A$
is not a $\pi$--base of $\U$. 
Without loss, all $A_\alpha$
are either contained in one $X_n$ or intersect each $X_n$ at most once.
For the second kind of $A_\alpha$, let $B_\alpha =
\{ n ; \; A_\alpha \cap X_n \neq \em \}$. There is $C \in \V$
such that the $A_\alpha$ of the first kind are not a $\pi$--base
of $\V_n$ for $n \in C$; let $D_n \in \V_n$ witness this.
Since the $B_\alpha$ don't
form a $\pi$--base of $\V$,  choose $E \sub C$ witnessing this.
We now see easily that $A_\alpha \not\sub^* \bigcup_{n \in E} \{ n \}
\times D_n 
\in \U$, as required.

(b) By (a) it suffices to prove that $\pi\chi (\U) \geq \pi\chi
(\W)$.  Let $\kappa < \pi\chi (\W)$,
and $\A=\{ A_\alpha ; \; \alpha < \kappa \} \sub [\omega\times\omega]^\omega
$. We want to show $\A$ is not a $\pi$--base of $\U$. 
For this simply let $B_\alpha  = \{ m ; \; \la n,m\ra \in A_\alpha$
for some $n\}$, find $C \in \W$ which does not almost contain
any of the $B_\alpha$ which are infinite, and note that $A_\alpha 
\not\sub^* \{ \la n , m \ra ; \; m>n$ and $m\in C \}\in \U$
(for any $\alpha$), as required. (Note that the same argument shows that
$\pi\chi (\U) \geq \pi\chi (f (\U))$ for all ultrafilters $\U$ and
all finite--to--one functions $f$.)

(c) This is easy, for given a $\pi\sigma$--base
$\A$ of $\U$, the family $\{ B_A ; \; A\in \A\}\cap\omoms$ where
$B_A = \{ n ; \; A\cap X_n \neq \em \}$ is a $\pi\sigma$--base
of $\V$.

(d)  By (c) we know
$\pi\chi_\sigma (\U) \geq \pi\chi_\sigma (\V)$;
$\pi\chi_\sigma (\U) \geq \pi\chi_\sigma (\W)$ is proved as
in (b); finally, $\pi\chi_\sigma (\U) \geq \dd$ follows
from Proposition 1.6 because $\U$ is not a $P$--point.
So we are left with showing that $\pi\chi_\sigma (\U)
\leq \max$.

For this choose
$\pi\sigma$--bases $\{ A_{\alpha} ; \; \alpha < \max \}$
of $\W$ and $\{ B_\beta ; \; \beta < \max \}$ of $\V$,
as well as a dominating family $\{ f_\gamma ; \; \gamma < \max\}$.
Put $C_{\alpha,\beta,\gamma} = \{ \la n,m \ra ; \;
n \in B_\beta$ and $m = \min (A_\alpha \sem f_\gamma (n) \}$.
To see that the $C_{\alpha, \beta , \gamma}$ form a $\pi\sigma$--base
of $\U$, take $D_k \in \U$, put $E_k = \{ n ; \; \{ m ;\;
\la n,m \ra \in D_k \} \in \W \} \in \V$, and find $\beta$
with $B_\beta \sub^* E_k$ for all $k$. Also let
$F_{k,n} = \{ m ;\;
\la n,m \ra \in D_k \}\in\W$ for $n \in E_k$, and find
$\alpha$ with $A_\alpha \sub^* F_{k,n}$ for all $k,n$.
Define $g_{k} \in \omom$ such that $A_\alpha \sub F_{k,n} \sem
g_k (n)$ for all $n$, and find $\gamma$ such that
$f_\gamma \geq^* g_k$ for all $k$. It is now
easy to see that $C_{\alpha,\beta,\gamma} \sub^*
D_k$ for all $k$. $\qed$

\bigskip

Another ultrafilter construction goes as follows:
let $\lambda$ be a regular uncountable cardinal,
let $\T = \la T_\alpha ; \; \alpha < \lambda \ra$
be a tower,
let $\V_\alpha$, $\alpha < \lambda$, be ultrafilters on
$\omega$ with $T_\alpha \in \V_\alpha$, and let
$\D$ be a uniform ultrafilter on $\lambda$. Define
$\U$ as the {\it $\D$--limit} of the $\V_\alpha$, i.e.
\sm

\ce{$U \in \U \Loleriar \{ \alpha ; \; U \in \V_\alpha \} \in \D$.}
\sm

\no Then:
\sm

{\capit Proposition 5.2.} {\it $\lambda \leq \pi\chi (\U) \leq \lambda
\cdot
\sum_\D \pi\chi (\V_\alpha)$.}
\sm

{\it Proof.} Note that $\T \sub \U$. Since $\T$ has no pseudointersection,
$\pi\chi(\U) \geq\lambda$ is immediate. To see the other inequality,
let $\A_\alpha $ be  $\pi$--bases of the $\V_\alpha$ for appropriate
$\alpha$'s. Then $\bigcup_\alpha \A_\alpha$ is a $\pi$--base
of $\U$. This shows $\pi\chi(\U) \leq \lambda
\cdot \sum_\D \pi\chi (\V_\alpha)$. $\qed$ 
\sm

\no As an immediate consequence we see
\sm

{\capit Corollary 5.3.} {\it Let $\kappa < \lambda$ be regular
uncountable cardinals such that there
is an ultrafilter $\V$ with $\pi\chi (\V) = \kappa$ and a 
tower of height $\lambda$. Then there is an ultrafilter
$\U$ with $\pi\chi (\U) = \lambda$.} $\qed$
\bigskip

{\bolds Theorem 4.} {\it (a) Let $R$ be a set of regular uncountable 
cardinals in $V\models GCH$. Then there is a forcing notion $\PP$
such that 
\sm

\ce{$V^\PP \models ``$for all $\lambda\in R$ there is an ultrafilter $\U$
such that $\pi\chi(\U) = \lambda$".}
\sm

(b) It is consistent there is an ultrafilter $\U$ with
$\pi\chi (\U) < \pi\chi_\sigma (\U)$. More explicitly,
given $\kappa < \lambda$ regular uncountable,
we can force $\pi\chi (\U) = \kappa$ and $\pi\chi_\sigma
(\U) = \lambda$ for some ultrafilter $\U$.

(c) It is consistent there is an ultrafilter $\U$
with $\pi\chi(\U)= $ \non$(r^0_\U) = \omega_\omega$. In particular
$\pi\chi(\U)=$ \non$(r^0_\U)$ consistently has countable cofinality.}
\sm

{\it Proof.} (a) We plan to adjoin an ultrafilter $\U$ with
$\pi\chi (\U) = \omega_1$ and towers $\T_\lambda$ of height $\lambda$
for each $\lambda \in R$. Then the result will follow by 5.3.
Note that the consistency of the existence of towers
of different heights was proved by Dordal [Do 1, section 2]
with essentially the same argument.

Let $\mu > \sup (R)^+$ be a regular
cardinal. We shall have $\PP = \PP^0 \times \PP^1$
where $\PP^0$ is the Easton product which adds $\mu$
subsets to $\lambda$ for each $\lambda \in R$ and $\PP^1$
is a $ccc$ forcing notion. Since $\PP^1$ is still $ccc$
in $V^{\PP^0}$, cofinalities and cardinals are preserved.

$\PP^1$ is an iteration $\PP^2 \star \dot\PP^3$
where $\PP^2$ is the finite support product of the
forcings $\QQ^\lambda$, $\lambda \in R$, which add
families $\{ C^\lambda_\eta ; \; \eta \in 2^{<\lambda} \}$
of infinite subsets of $\omega$
such that
\sm

\item{$\bullet$} $\eta \sub \theta$ implies $C^\lambda_\theta
\sub^* C^\lambda_\eta$ and

\item{$\bullet$} $C^\lambda_{\eta \ha\la 0\ra} \cap
C^\lambda_{\eta\ha\la 1 \ra}$ is finite,
\sm

\no with finite conditions (see the proof of Theorem 5 for
a similar, but more complicated, forcing).
In $V^{\PP^2}$, $\PP^3$ is a finite support iteration of
length $\omega_1$ of Mathias forcings with an ultrafilter
(see $\S\S$ 2 and 4) which adds an ultrafilter $\U$
all of whose cardinal coefficients are equal to $\omega_1$.
(Alternatively, we could define $\PP^3$ in $V^{\PP^0 \times
\PP^2}$.)

Since $\PP^0$ is still $\omega_1$--distributive over $V^{\PP^1}$,
it doesn't add reals, and $\U$ still is an ultrafilter
with $\pi\chi(\U) = \omega_1$ in $V^{\PP}$ 
($\star$). Also, 
in $V^\PP$ we have $\cc < \mu$, but
$2^\lambda = \mu$ for all $\lambda \in R$.
Given $f \in 2^\lambda$ with $f\re\alpha \in V$ for
all $\alpha < \lambda$, 
\sm

\ce{$\T_f = \{ C^\lambda_{f\re \alpha} ; \; \alpha < \lambda \}$}
\sm

\no forms a $\sub^*$--decreasing chain. Because of the Easton product,
we have $\mu$ such $\T_f$'s for each $\lambda$. Since $\cc < \mu$,
not all of them can have a pseudointersection. Hence,
for each $\lambda \in R$, there is a tower $\T_\lambda$ of height $\lambda$ 
($\star
\star$). (In fact, a density argument shows none of them
has a pseudointersection, see the proof of Theorem 5.)
By Corollary 5.3 as well as $(\star)$ and $(\star\star)$,
we have, for each $\lambda \in R$, an ultrafilter $\V$
with $\pi\chi (\V) = \lambda$.
\sm

(b) 
By [BlS] we know it is consistent there is a $P$--point
$\V$  with $\kappa =\pi\chi (\V) =\pi\chi_\sigma (\V) < \dd =\lambda$. 
Put $\U = \V \times \V$.
Then $\pi\chi (\U) = \pi\chi(\V) = \kappa$ by 5.1 (b)
and $\pi\chi_\sigma (\U) = \dd = \lambda$ by 5.1 (d).
(Instead of [BlS], the ultrafilters gotten in the construction
in Theorem 5 could be used for this consistency.) 
\sm

(c) By (a) we can force ultrafilters $\V_\alpha$ for all regular
$\alpha$ with $\omega_1 \leq\alpha \leq \omega_{\omega+1}$.
Then $\U = \sum_{\V_{\omega_{\omega+1}}} \V_n$ satisfies
$\pi\chi (\U) = \omega_{\omega}$, by 5.1 (a). $\qed$
\bigskip

We conclude this section with the discussion of several
refinements of Theorem 4.
The construction in part (b) of the above
proof also shows that the result in Lemma 1.6
is sharp and cannot be improved.
\sm

{\capit Corollary 5.4.} {\it It is consistent there is
an ultrafilter $\U$ which is not a $P$--point such that
$\pi\chi (\U) < \dd$.} $\qed$
\sm

\no The result in Theorem 4 (a) will be superseded by Theorem 5 in
the next section. We still gave its proof because it is
much simpler and also because of the following two consequences
of the construction which we cannot get from Theorem 5.

\sm

{\capit Corollary 5.5.} {\it In the statement of Theorem 4 (a), we can
delete the word ``regular".}
\sm

{\it Proof.} Assume without loss that whenever $\lambda \in R$
is singular, then $R \cap \lambda$ is cofinal in $\lambda$.
We show the construction in the proof of Theorem 4 (a) produces
an ultrafilter $\U$ with $\pi\chi (\U) = \lambda$.

For $\mu \in R \cap \lambda$ we added towers $\T_\mu = \la
T_{\alpha , \mu} ; \; \alpha < \mu\ra$ such that
$\bigcup_{\mu \in R \cap \lambda} \T_\mu$ is a filter base (this is 
immediate from the definition of the forcing $\PP^2$).
Put $S = \{ F ; \; F$ is a finite subset of $\{ \la \alpha , \mu \ra;
\; \alpha < \mu$ and $\mu \in R \cap \lambda \}\}$.
Clearly $|S| = \lambda$. For $F \in S$ let $\V_F$ be an ultrafilter
on $\bigcap_{\la \alpha , \mu\ra \in F} T_{\alpha , \mu}$
with $\pi\chi (\V_F) = \omega_1$. Let $\D$ be a uniform ultrafilter on
$S$ such that for any $F \in S$, $\{ G \in S ; \; G \supseteq
F \} \in \D$. Put $\U = \lim_\D \V_F = \{ X \sub \omega ;
\; \{ F \in S ; \; X \in \V_F \} \in \D \}$. We have to show
$\pi\chi (\U) = \lambda$. This is no more but an elaboration of the
argument in Proposition 5.2.

$\pi\chi (\U) \leq\lambda$ is immediate since the union of the 
$\pi$--bases of the $\V_F$ is a $\pi$--base of $\U$. To see
$\pi\chi (\U) \geq \lambda$, it suffices to show that
$\T_\mu \sub \U$ for each $\mu \in R\cap \lambda$.
Fix $\alpha < \mu$. As $\{ G \in S ; \; G \ni
\la \alpha , \mu \ra \} \in \D$ and $T_{\alpha,\mu} \in \V_G$
for any $G\in S$ with $\la \alpha , \mu\ra \in G$, it follows that
$T_{\alpha , \mu} \in \U$, as required. $\qed$
\sm

{\capit Corollary 5.6.} {\it In Theorem 4 (a), we can additionally
demand that the dominating number $\dd$ is an arbitrary
regular uncountable cardinal. In particular, there may be
many different $\pi$--characters below $\dd$.}
\sm

{\it Proof.} Simply replace the forcing $\PP^3$ in the proof of Theorem 4 (a)
by the forcing from [BlS] which adds an ultrafilter $\U$ with
$\chi (\U) = \omega_1$ while forcing $\dd = \kappa$ where
$\kappa$ is an arbitrary regular cardinal. $\qed$

\Bigskip


{\dunhg 6. The spectral problem}
\Smallskip

\no By Louveau's Theorem mentioned in $\S$ 1, we know it is consistent
that there are simultaneously ultrafilters with many different
values for $\pp$. The same is true for $\pi\chi$, as proved in
the preceding section.
We now complete this cycle of results 
by showing how to get the 
consistency of the simultaneous existence of many ultrafilter characters
and, dually, of many values for $\pi\pp$.
\bigskip

{\bolds Theorem 5.} {\it Let $R$ be a set of regular uncountable
cardinals in $V \models GCH$. Then there is a forcing notion
$\PP$ such that
\sm

\ce{$V^\PP \models$ ``for all $\lambda \in R$ there is 
an ultrafilter $\U$ such that
$\chi (\U) = \pi\chi(\U) =\lambda$".}}

\sm

In fact, the ultrafilters we construct in the proof are all $P$--points.

\sm

{\it Proof.} We plan to adjoin, for each $\lambda \in R$,
a matrix $\la E^\alpha_{\lambda,\gamma} ; \; \alpha < \omega_1
, \gamma < \lambda \ra$ of subsets of $\omega$ such that
the following conditions are met:
\sm

\item{(i)} $\la E^\alpha_{\lambda,\gamma} ; \; \gamma < \lambda \ra$
forms a tower;

\item{(ii)} $\alpha < \beta < \omega_1$ entails $E^\beta_{\lambda,
\gamma} \subseteq^* E^\alpha_{\lambda,\gamma}$;

\item{(iii)} for each $X \sub \omega$ we find a pair $\la \alpha,
\gamma\ra$ such that either $E^\alpha_{\lambda,\gamma} \subseteq^*
X$ or $E^\alpha_{\lambda,\gamma} \subseteq^* \omega\sem X$.
\sm

\no Clearly, this is enough: all three conditions imply 
the matrix generates an ultrafilter, we get $\chi(\U) 
\leq \lambda$ by the size of the matrix, and (i) entails
$\pi\chi (\U) \geq\lambda$.

We now describe the forcing we use. We shall have
$\PP = \PP^0 \times \PP^1$ where $\PP^1 $ is $ccc$ and
$\PP^0$ is the Easton product of the forcing notions adding
one subset of $\lambda$ with conditions of size $<\lambda$
for $\lambda \in R$.
Since $\PP^0$ is $\omega_1$--closed, it preserves the $ccc$
of $\PP^1$, and thus $\PP$ preserves cofinalities and cardinals.
However, we shall look at $\PP$ as first forcing with $\PP^1$
and then with $\PP^0$. In $V^{\PP^1}$, the closure 
property of $\PP^0$
is lost, but it is still $\omega_1$--distributive.

To define $\PP^1$, put $$\mu = \cases{ \sup (R) & if this has
uncountable cofinality \cr
\sup (R)^+ & otherwise \cr}$$
Let $\zeta = \mu \cdot \omega_1$, and let $\{ A_\lambda ; \;
\lambda \in R\}$ be a partition of $\zeta$ such that
$|A_\lambda \cap [\mu\cdot\beta , \mu (\beta + 1))| = \mu$ for each $\beta
< \omega_1$.
$\PP^1$ shall be a finite support iteration 
$\la \PP_\alpha , \dot \QQ_\alpha ; \; \alpha < \zeta \ra$ of
$ccc$ p.o.'s such that
\sm

\ce{$\forces_{\alpha} ``| \dot \QQ_\alpha |\leq \mu"$}
\sm

\no for all $\alpha < \zeta$. Since we have $GCH$ in the ground
model, this implies $V^{\PP^1} \models \cc \leq \mu$
so that we can enumerate the names of subsets of $\omega$
arising in the extension in order type $\zeta$.
More explicitly, we shall have a sequence $\la \dot X_\alpha ; 
\; \alpha < \zeta \ra$ such that
\sm

\item{$\bullet$} $\forces_\alpha ``\dot X_\alpha \sub \omega"$ and

\item{$\bullet$} whenever $\dot X$ is a $\PP_\alpha$--name
for a subset of $\omega$ and $\lambda \in R$, then there is $
\beta \geq \alpha$, $\beta \in A_\lambda$, such that
\sm

\ce{$\forces_\beta `` \dot X = \dot X_\beta"$.}
\sm

\no Clearly this can be done.

Along the iteration, we want to add, for $\lambda \in R$ and
$\alpha \in A_\lambda$, a system
$\la C^\alpha_\eta ; \; \eta \in 2^{< \lambda } \cap V \ra$
of infinite subsets of $\omega$ lying in $V^{\PP_{\alpha + 1}}$ such that
\sm

\item{(1)} $\eta \sub \theta$ implies $C^\alpha_\theta \sub^*
C^\alpha_\eta$;

\item{(2)} $C^\alpha_{\eta \ha\la 0\ra} \cap C^\alpha_{\eta \ha
\la 1 \ra}$ is finite;

\item{(3)} if $\alpha < \beta$ both belong to $A_\lambda$,
then $C^\beta_\eta \sub^* C^\alpha_\eta$;

\item{(4)} in $V^{\PP_{\alpha + 1}}$, the set
$\{ \eta \in 2^{< \lambda} \cap V ; \; C^\alpha_\eta \sub^* X_\alpha$
or $C^\alpha_\eta \cap X_\alpha$ is finite$\}$ is dense
(and open) in $2^{<\lambda} \cap V$.
\sm

\no In (4), $X_\alpha$ denotes, of course, the interpretation
of $\dot X_\alpha$ in $V^{\PP_{\alpha+1}}$.

We are ready to describe the factors $\QQ_\alpha$ of the iteration.
Fix $\alpha$ and work in $V^{\PP_\alpha}$. We distinguish
two cases:
\sm

{\sanse Case 1.} $\alpha = \min (A_\lambda)$ or $cf (A_\lambda
\cap \alpha) \geq \omega$. 
\sm

\no $\QQ_\alpha$ consists of pairs $\la s , F \ra$ where
$F \sub A_\lambda \cap \alpha$ is finite (the second
part of the condition is missing in case
$\alpha = \min (A_\lambda)$) and  $s$ is a finite partial
function with $dom (s) \sub 2^{<\lambda } \cap V$ and such that
$s (\eta) \sub \omega$ is finite for all $\eta \in dom (s)$. 
We stipulate $\la s, F \ra \leq \la t , G \ra$ iff
$G \sub F$, $dom (t) \sub dom (s)$ and
$t(\eta) \sub s(\eta)$ for all $\eta\in dom (t)$ as well as \sm

\item{$\bullet_1$} if $\eta \sub \theta$ belong to $ dom (t)$,
then $s(\theta) \sem t(\theta) \sub s(\eta) \sem t(\eta)$;

\item{$\bullet_2$} if $\eta\ha\la 0 \ra , \eta\ha\la 1 \ra
\in dom (t)$, then $s(\eta\ha\la 0 \ra) \sem t(\eta\ha\la 0 \ra)$
and $s(\eta\ha\la 1 \ra) \sem t(\eta\ha\la 1 \ra)$ are disjoint;

\item{$\bullet_3$} if $\alpha \in G$ and $\eta \in dom (t)$, then
$s(\eta) \sem t(\eta) \sub C^\alpha_\eta$.
\sm

\no This forcing is easily seen to be $ccc$, and it adds
a system $\la D^\alpha_\eta ; \; \eta \in 2^{<\lambda } \cap
V\ra$ of subsets of $\omega$ which satisfies (1) through (3)
above by $\bullet_1$ through $\bullet_3$.
\sm

{\sanse Case 2.} $\beta = \max (\alpha \cap A_\lambda)$.
\sm

\no Let $\QQ_\alpha$ be the trivial ordering, and define
$D^\alpha_\eta : = C^\beta_\eta$ for all $\eta \in 2^{<\lambda}
\cap V$.

This completes the construction of the $\QQ_\alpha$, and, hence,
of the forcing $\PP^1$. We still have to explain how
to get the $C^\alpha_\eta$ from the $D^\alpha_\eta$ in the
model $V^{\PP_{\alpha + 1}}$. For this, note that the
set $$\eqalign{\{ \eta \in 2^{< \lambda } \cap V ; \; &
D^\alpha_\eta \sub^* X_\alpha \hbox{ or } D^\alpha_\eta \cap X_\alpha
\hbox{ is finite or} \cr
& \hbox{ for all } \theta \supseteq \eta, \hbox{ both }
D^\alpha_\theta \cap X_\alpha \hbox{ and } D^\alpha_\theta \cap
(\omega \sem X_\alpha ) \hbox{ are infinite}\} \cr }$$
is dense and open in $2^{< \lambda} \cap V$.
Let
$$C^\alpha_\eta = \cases{ D^\alpha_\eta \cap X_\alpha &
if $\eta$ enjoys the third property \cr
D^\alpha_\eta & otherwise \cr}$$
Then (2) and (3) are trivially true, and it is easy to check
that (1) and (4) are satisfied as well. Thus we are done with
the construction of the required system.

Next, let $f_\lambda \in 2^\lambda$, $\lambda \in R$,
be the generic Easton functions. Also let
$B_\lambda \sub A_\lambda$ be a cofinal subset of order type
$\omega_1$. For $\alpha < \omega_1$, set
$E^\alpha_{\lambda , \gamma} = C^{B_\lambda (\alpha)}_{f_\lambda
\re \gamma}$ where $B_\lambda (\alpha)$ denotes the $\alpha$--th
element of $B_\lambda$. We claim the $E^\alpha_{\lambda, \gamma}$
satisfy (i) through (iii) above.

Now, (ii) is immediate from (3). To see (i), first note that
the $E^\alpha_{\lambda , \gamma}$ for fixed $\alpha$
form a decreasing chain, by (1). Next use
a genericity argument to see that this
chain has no pseudointersection, as follows. 
Work in $V^{\PP^1}$. By distributivity (see above),
$\PP^0$ adds no new reals over $V^{\PP^1}$. 
Hence it suffices to check that given any $X \sub \omega$ in
$V^{\PP^1}$, the set $\{ \eta \in 2^{< \lambda} ; \;
X \not\sub^* C^\alpha_\eta \}$ is dense in $2^{< \lambda}
\cap V$. But this is trivial by (2). Finally,
(iii) is taken care of by an exactly analogous density argument
involving (4). Hence we're done. $\qed$
\bigskip

In the remainder of this section,
we discuss several improvements of, and variations on, 
the above result which
are corollaries to the construction.
\sm

{\capit Corollary 6.1.} {\it If we care only about characters,
we can relax the assumption about $R$ in Theorem 5 to:
``$R$ is a set of cardinals of uncountable cofinality."}

\sm

{\it Proof.} We confine ourselves to describing the changes we need
to make in the proof of Theorem 5. We additionally adjoin, for 
$\lambda \in R$ singular, a matrix $\la E^\alpha_{\lambda , \Gamma}
; \; \alpha < \omega_1 , \Gamma \in [\lambda]^\omega \cap V \ra$
of subsets of $\omega$ such  that, in addition to (ii) and (iii)
(with $\gamma $ replaced by $\Gamma$), the following conditions
are met:
\sm

\item{(i$_a$)} $\Gamma \sub \Delta$ implies $E^\alpha_{\lambda ,
\Delta} \sub^* E^\alpha_{\lambda , \Gamma}$;

\item{(i$_b$)} $\Gamma \not\sub\Delta$ implies $E^\alpha_{\lambda ,
\Delta} \not\sub^* E^\beta_{\lambda , \Gamma}$ for all $\alpha
, \beta$.
\sm

\no Then (i$_a$), (ii) and (iii) imply the matrix generates an
ultrafilter $\U$ with $\chi (\U) \leq \lambda$, while (i$_b$) gives us 
that $\chi (\U) \geq \lambda$.

With $\PP^0$ we also adjoin a function $f_\lambda$ from $\lambda$
to $2$ with countable conditions for each singular $\lambda \in R$.
In the partition of $\zeta$ include the $A_\lambda$ for singular
$\lambda$. $\PP^1$ is as before except that we still have
to define $\QQ_\alpha$ for $\alpha \in A_\lambda$
where $\lambda$ is singular.

For such $\lambda$ and $\alpha \in A_\lambda$, we add a system
$\la C^\alpha_\eta ; \; \eta \in V , \eta : \lambda \to 2 $
is a partial function with countable domain$\ra$ of subsets of
$\omega$ such that, in addition to (3) and (4), we have
\sm

\item{(1$_a$)} $\eta \sub \theta$ implies $C^\alpha_\theta \sub^* 
C^\alpha_\eta$;

\item{(1$_b$)} $\eta \not\sub \theta$ implies $C^\alpha_\theta \not\sub^*
C^\beta_\eta$ for all $\alpha , \beta$;

\item{(2)} if $\eta$ and $\theta$ are incompatible, then
$C^\alpha_\eta \cap C^\alpha_\theta$ is finite.
\sm

\no The corresponding $D^\alpha_\eta$ are produced as before
and satisfy (1$_a$), (1$_b$), (2) and (3). Note that
(1$_b$) for the $D^\alpha_\eta$ is easily preserved 
in Case 1 by a genericity argument. $C^\alpha_\eta$ is defined
from $D^\alpha_\eta$ as previously and satisfies
(1$_a$), (2), (3) and (4). To see that is also satisfies
(1$_b$) suppose that $C^\alpha_\theta \sub^* C^\beta_\eta$
for some $\alpha \geq\beta$ and $\eta \not\sub \theta$. Then
$\theta$ has an extension $\bar\theta$ such that
$\eta$ and $\bar\theta$ are incompatible and
thus $| C_\eta^\beta \cap C_{\bar\theta}^\beta | < \omega$
by (2). By (1$_a$), $C^\alpha_{\bar\theta} \sub^* C_\theta^\alpha$
which means that $| C^\alpha_{\bar\theta} \cap C_{\bar\theta}^\beta
| < \omega$, contradicting (3). (This is the only place where
(2) is needed.)

Finally put $E^\alpha_{\lambda, \Gamma } = C^{B_\lambda (\alpha)}_{
f_\lambda \re \Gamma }$ and use (1$_a$), (1$_b$), (3) and (4)
to check that (i$_a$), (i$_b$), (ii) and (iii)
are satisfied. $\qed$

\sm

{\capit Corollary 6.2.} {\it In Theorem 5, we can additionally demand
that all the ultrafilters produced are Ramsey ultrafilters.}
\sm

{\it Proof.} We replace the $\QQ_\alpha$ in the iteration by
$\QQ_\alpha \star \CC$ where $\CC$ denotes Cohen forcing
(so $\PP_{\alpha+1} = \PP_\alpha \star\dot \QQ_\alpha \star\CC$).
Apart from that, it suffices to change the way the
$C^\alpha_\eta$ are defined from the $D^\alpha_\eta$.
Instead of listing names for subsets of $\omega$, we list
names for partitions of $\omega$, as $\la \dot X_{\alpha , n} ; 
\; n\in\omega \ra$. Assume we are at step $\alpha$. Look at
$$\eqalign{\{ \eta \in 2^{< \lambda } \cap V ; \; &
D^\alpha_\eta \sub^* \bigcup_{n<N} X_{\alpha ,n} \hbox{ for some }
N \hbox{ or }  \cr
& D^\alpha_\theta \hbox{ meets infinitely many } X_{\alpha
, n} \hbox{ for all } \theta \supseteq \eta\} \cr }$$
This is again dense and open in $2^{< \lambda} \cap V$.
Now let $Y$ be a Cohen real over $V^{\PP_{\alpha} \star \dot\QQ_\alpha}$,
and think of $Y$ as a subset of $\omega$ which meets
each $X_{\alpha , n}$ once. Then let
$$C^\alpha_\eta = \cases{ D^\alpha_\eta \cap Y &
if $\eta$ enjoys the second property above\cr
D^\alpha_\eta & otherwise \cr}$$
We now see that the $C^\alpha_\eta$ satisfy
\sm

\item{(4')} in $V^{\PP_{\alpha + 1}}$, the set
$\{ \eta \in 2^{< \lambda} \cap V ; \; C^\alpha_\eta \sub^* 
\bigcup_{n<N} X_{\alpha,n}$ for some $N$
or $|C^\alpha_\eta \cap X_{\alpha,n}|\leq 1$ for all $n\}$ is dense
(and open) in $2^{<\lambda} \cap V$.
\sm

\no And thus we get by genericity
\sm

\item{(iii')} for each partition $\la X_n ; \; n\in \omega\ra$ of $\omega$,
we find a pair $\la \alpha,
\gamma\ra$ such that either $E^\alpha_{\lambda,\gamma} \subseteq^*
\bigcup_{n<N}
X_n$ for some $N$ or $|E^\alpha_{\lambda,\gamma} \cap X_n|
\leq 1$ for all $n$,
\sm

\no which guarantees Ramseyness. $\qed$
\bigskip

\ce{$\star\star\star$}
\bigskip

\no By Theorem 5, we can get a plethora of ultrafilter
characters. This suggests it might be interesting to know
whether an arbitrary set of regular cardinals can be realized
as the set of possible ultrafilter characters in some
model of $ZFC$. To this end
we define
\sm

\item{$\bullet$} \Spec {\chi} $= \{ \lambda ; \; \chi (\U) = \lambda$
for some ultrafilter $\U$ on $\omega\}$, the {\it character
spectrum};

\item{$\bullet$} \Spec {\pi\chi} $= \{ \lambda ; \; \pi\chi (\U) = \lambda$
for some ultrafilter $\U$ on $\omega\}$, the {\it $\pi$--character
spectrum};

\item{$\bullet$} \Specs {\chi} $=\{ \lambda ; \; \chi (\U) 
=\pi\chi (\U) = \lambda$
for some ultrafilter $\U$ on $\omega\}$.
\sm

\no Unfortunately, we have no limitative results on \Spec {\chi} and
on \Spec {\pi\chi} (see section 8 for some questions on this;
in particular, question (5)),
but we can prove the following which answers the spectral
question for \Specs {\chi} in many cases.
\bigskip

{\bolds Theorem 6.} {\it Let $R$ be a non--empty set of uncountable
regular cardinals in $V\models GCH$. Then there is a forcing notion
$\PP$ such that in $V^\PP$, for all regular $\lambda$ which are neither
successors of singular limits of $R$ nor inaccessible limits of $R$,
we have
\sm

\ce{$\lambda \in R \Loleriar \lambda \in$ \Specs \chi.}}
\sm

{\it Proof.} We use a modification of the partial order in Theorem 5,
and confine ourselves to describing the differences of the two proofs.
Let $\nu = \min (R)$. Put $\zeta = \mu \cdot \nu$.
$\PP^0$ is defined exactly as before, and $\PP^1$ is
a finite support iteration of length $\zeta$ which
\sm

\item{(a)} takes care of all the p.o.'s described in the proof of
Theorem 5; and

\item{(b)} forces $MA$ for all p.o.'s of size $<\mu$ at each limit
step of the form $\mu \cdot \beta$ where $\beta < \nu$ is a successor
ordinal.
\sm

\no It's clear that this can be done. By Theorem 5, we know that
$R \sub$ \Specs \chi. We proceed to show the other direction.

Let $\lambda\not\in R$ be a regular cardinal which is neither a
successor of a singular limit of $R$ nor an inaccessible limit
of $R$. Assume there is, in $V^\PP$, an ultrafilter $\U$
with $\chi (\U) = \lambda$. Since $\cc = \mu$, and
$\mu$ doesn't qualify as $\lambda$ (because either $\mu$
is a limit of $R$ (and thus either inaccessible or not regular)
or $\mu$ is a successor of a singular limit of $R$ or
$\mu = \sup (R) = \max (R) \in R$), we know $\lambda < \mu$.
Since the cofinality of the iteration is $\nu \in R$,
we also see $\lambda > \nu$.
We shall show that $\pi\chi(\U) < \lambda$.
Let $\F \sub \U$ be a base of $\U$ of size $\lambda$.
Work in $V^{\PP^1}$. The forcing $\PP^0$ decomposes as
a product $\PP^{<\lambda} \times \PP^{>\lambda}$ because $\lambda
\not\in R$. The first part has size $<\lambda$. This follows from
the Easton support in case $\lambda$ is a successor of an inaccessible,
and is trivial in the other cases. The second part is
$ \lambda^+$--distributive, and thus adds no new sets
of size $\lambda$. Hence $\F \in V^{\PP^{<\lambda} \times
\PP^1}$. 

For each $p \in \PP^{<\lambda}$, let (in $V^{\PP^1}$)
$\F_p = \{ F \sub \omega ; \; p \forces_{\PP^{>\lambda}} F \in \dot\F \}$.
Clearly, $\F \sub \bigcup_p \F_p$.
Next, for each $\beta < \nu$, let (in $V^{\PP_{\mu\cdot\beta}}$)
$\F_{p,\beta} = \{ F \sub \omega ; \; F \in V^{\PP_{\mu\cdot\beta}}$
and $\forces_{\PP^1 / \PP_{\mu\cdot\beta}} F \in \dot\F_p \}$.
By $ccc$--ness of the iteration, we have $\F_p = \bigcup_\beta
\F_{p,\beta}$. Since we forced $MA$ along the way, each
$\F_{p,\beta}$ has a pseudointersection which we call $G_{p, \beta}$.
By construction, the $G_{p , \beta}$ form a $\pi$--base
of $\U$ which has size $< \lambda$, as required. $\qed$
\sm

\no Note the similarity between this result and Dordal's
result [Do 1, Corollary 2.6] on the spectrum of tower heights.
The latter is easier to prove because a tower is an easier
combinatorial object than an ultrafilter. 
---
Of course, the restrictions on $\lambda$ in the above theorem
come from the size of the set of Easton conditions, and the present
proof does not work in the other cases. In the other direction,
we can show that certain cardinals may arise as characters
even if they don't belong to $R$:

\sm

{\capit Proposition 6.3.} {\it In the model constructed in Theorem 6: if
$R$ contains cofinally many $\omega_n$, then $\omega_{\omega + 1}
\in $ \Spec \chi.}
\sm

{\it Proof.} Let $S = R \cap \omega_\omega$,
let $\nu = \min (R)$ as before, and put  $S' = \{ n ; \;
\omega_n \in S\}$. Assume
we have, for $n \in S'$, an ultrafilter $\V_n$ on $\omega$
which is generated by a matrix
$\la E^\alpha_{\omega_n,\gamma} ; \; \alpha < \nu
, \gamma < \omega_n \ra$ satisfying (i) through (iii) in the proof
of Theorem 5. Put $\U : = \sum_\V \V_n$ where
$\V$ is any ultrafilter with $\chi (\V) = \nu < \omega_\omega$
and $S' \in \V$. We shall show that $\chi (\U) = \omega_{\omega+1}$.

To see $\chi (\U) \geq \omega_{\omega+1}$, it suffices to show that
$\chi (\U) \geq \omega_n$ for all $n$, by Proposition 1.4.
This is easy: fix $n $ and $\F \sub \U$ with $|\F| = \omega_n$;
for $m > n$ with $m \in S'$ find $A_m \in \V_m$
such that $F \cap ( \{ m \} \times \omega) \not\sub^* \{ m\}
\times A_m$ for any $F \in \F$ with $F \cap ( \{ m \} \times \omega )
\in \V_m$, and put $A = \bigcup_{m > n , m \in S'}
\{ m \} \times A_m \in \U$; then $F \not\sub^* A$ for any
$F \in \F$, as required.

To see $\chi (\U) \leq \omega_{\omega + 1}$, note first that
$\dd = \nu$ by construction. Now, let $\{ g_\delta ; \; \delta < \nu \}$
be a dominating family, and let $\{ V_\zeta ; \; \zeta < \nu \}$ be
a base of $\V$; without loss, each $V \in \V$  strictly
contains at least one $V_\zeta$. Also let
$\{ f_\eta ; \; \eta < \omega_{\omega+1} \} \sub \prod S: =
\{ f: S ' \to \omega_\omega ; \; f(n) < \omega_n \}$ be a dominating
family; i.e. given $f \in \prod S$ there is $\eta < \omega_{\omega +1}$
such that $f(n) < f_\eta (n)$ for all $n \in S'$
(such a family clearly exists in the ground model; it also
exists in the generic extension because $pcf$ is left unchanged
by the forcing).
For $\alpha , \delta , \zeta < \nu$ and $\eta < \omega_{
\omega + 1}$, put $A_{\alpha , \delta , \zeta , \eta } = \bigcup_{
n \in V_\zeta} \{ n \} \times ( E^\alpha_{\omega_n , f_\eta (n)}
\sem g_\delta (n))$ and check that the $A_{\alpha , \delta , \zeta,
\eta}$ form a base of $\U$. $\qed$
\sm

\no Note that for $\chi (\U) \geq \omega_{\omega + 1}$
we used no extra assumptions while the proof of $\chi (\U) \leq
\omega_{\omega + 1}$ involved the special shape of the ultrafilters $\V_n$
as well as $\dd < \omega_\omega$ (which is necessary by 1.6).
We don't know whether a similar result can be proved
without these assumptions (see $\S$ 8 (6)). 
Also, contrary to the situation for $\pi\chi$ (see Corollary 5.6),
we don't know whether we can have many characters below $\dd$.
\bigskip

\ce{$\star\star\star$}

\bigskip

\no We finally come to the result dual to Theorem 5.
\bigskip

{\bolds Theorem 7.} {\it Let $R$ be a set of regular
uncountable cardinals in $V \models GCH$. Then there is a
forcing notion $\PP$ such that
\sm

\ce{$V^\PP \models$ `` for all $\lambda \in R$ there is
an ultrafilter $\U$ such that $\pp (\U) = \pi\pp (\U) = \lambda$".}
}
\sm

{\it Proof.} Again, let
$$\mu=\cases{ \sup (R) & if this has uncountable cofinality \cr
\sup(R)^+ & otherwise \cr}$$
and adjoin, for all $\lambda \in R$, matrices
$\la E^\alpha_{\lambda , \gamma} ; \; \alpha < \mu , \gamma < \lambda \ra$
of subsets of $\omega$ such that (i) thru (iii) in the proof of
Theorem 5 are satisfied with $\omega_1$ replaced by $\mu$.
It is immediate that the matrices will generate ultrafilters of the required
sort. The rest of the proof of Theorem 5 carries over with very minor
changes which we leave to the reader. $\qed$
\sm

\no As in Corollary 6.2 we get
\sm

{\capit Corollary 6.4.} {\it In Theorem 7, we can additionally
demand that all the ultrafilters produced are Ramsey ultrafilters.}
$\qed$

\sm

One can again define \Spec \pp , \Spec{\pi\pp} and \Specs{\pp}
in the obvious fashion, but we do not know of any restrictive
results (like, e.g., Theorem 6) concerning these
spectra. ---
The second part of the following corollary --- which is immediate
from Theorems 5 and 7, Corollaries 6.2 and 6.4,
and results mentioned in $\S$ 1 --- answers the other half of
the question in [Br, subsection 4.1].
\sm

{\capit Corollary 6.5.} {\it (a) It is consistent with $ZFC$ that for some
Ramsey ultrafilter $\U$, \cov$(r^0_\U) = \pi\pp (\U) < \ppar$. \par
(b) It is consistent with $ZFC$ that for some Ramsey ultrafilter $\U$,
\non$(r^0_\U) = \pi\chi (\U) > \hhom$.} $\qed$

\Bigskip


{\dunhg 7. Connection with reaping and splitting}
\Smallskip

\no As we remarked in $\S$ 1, for any ultrafilter $\U$ on $\omega$
we have $\pi\pp (\U) \leq \ss$ and $\pi\chi (\U) \geq \rr$.
Furthermore, it follows from the results in $\S\S$ 5 and 6 that it
is consistent to have an ultrafilter $\U$ with $\pi\pp (\U) < \ss$,
as well as to have one with $\pi\chi (\U) > \rr$. Still there is
a close connection between the $\pi\chi (\U)$ and $\rr$,
as shown by the following well--known result whose proof we repeat
for completeness' sake.
\sm

{\capit Proposition 7.1.} (Balcar--Simon, [BS, Theorem 1.7])  $\rr = \min_\U
\pi\chi (\U)$.
\sm

{\it Proof.} Let $\A$ be a reaping family of size $\rr$. Without loss,
$\A$ is downward closed, that is, whenever $A \in \A$, then $\{ B \in
\A ; \; B \sub A \}$ is a reaping family inside $A$. This entails 
it can be shown by induction that given
pairwise disjoint $X_i$, $i \in n$, with $\bigcup_{i\in n} X_i =
\omega$, there are $i\in n$ and $A \in \A$ with $A \sub X_i$ $(\star)$. 
Let $\I$ be the ideal generated by sets $X \sub \omega$
with $A \not\sub^* X$ for all $A \in \A$. By $(\star)$, $\I$ is a proper
ideal. Hence it can be extended to a maximal
ideal whose dual ultrafilter $\U$ has $\A$ as a $\pi$--base, and thus
witnesses $\pi\chi (\U) = \rr$. $\qed$

\sm

Let us note that Balcar and Simon proved a much more general
result: the analogue of 7.1 holds in fact for a large class
of Boolean algebras.

We shall now see that there is no dual form of this proposition.
\bigskip

{\bolds Theorem 8.} {\it It is consistent with $ZFC$ that
$\pi\pp (\U) = \omega_1$ for all ultrafilters $\U$ on $\omega$,
yet $\ss = \cc = \omega_2$.}

\bigskip

For the proof of this Theorem we need to introduce several notions
and prove a few preliminary Lemmata.
Given a limit ordinal $\delta < \omega_2$, let $\la \delta_\zeta ; \;
\zeta \in cf(\delta) \ra$ be a fixed continuously increasing
sequence with $\delta = \bigcup_\zeta \delta_\zeta$. We define
sequences $\bar A^\alpha = \la A^\alpha_\beta \sub \alpha ; \; 
\beta < \omega_1 \ra$ for $\alpha < \omega_2$ recursively as follows.
\sm

\item{$\bullet$} $A^0_\beta = \em$

\item{$\bullet$} $A^{\alpha + 1}_\beta = A^\alpha_\beta \cup
\{ \alpha \}$

\item{$\bullet$} $A^\delta_\beta = \{ \gamma < \delta ; \; \gamma
\in A^{\delta_\zeta}_\beta$ for all $\zeta$ with $\delta_\zeta
> \gamma \}$ in case $cf(\delta) = \omega$

\item{$\bullet$} $A^\delta_\beta = \{ \gamma < \delta ; \;
\gamma < \delta_\zeta$ for some $\zeta < \beta$,
and $\gamma \in A^{\delta_\zeta}_\beta$
for all $\zeta < \beta$ with $\delta_\zeta > \gamma \}$
in case $cf(\delta) = \omega_1$
\sm

\no We leave it to the reader to verify that all $A^\alpha_\beta$
are at most countable and that for all
$\gamma < \alpha$, the set $\{\beta <\omega_1 ;
\; \gamma \in A_\beta^\alpha \}$ contains a club.

Now,
suppose $\FFF = \{ \F^\gamma = \{ F^\gamma_\beta ; \; \beta < \omega_1 \}
; \; \gamma < \alpha \}$ is a family of filter bases on $\omega$.
We call $\FFF$ {\it $\alpha$--nice} (or simply {\it nice}
if the $\alpha$ in question is clear from the context)
iff, given $X \in \omoms$
and a set $\{ f_j ; \; j\in\omega\}$ of one--to--one
functions in $\omom$, there
is a club $C = C(X , \la f_j \ra_j) \sub \omega_1$ such that $| X \sem \bigcup_{j < k}
\bigcup_{\gamma \in \Gamma}
f_j^{-1} (F_\beta^\gamma) | = \omega$ for all $\beta \in C$, 
all $k\in\omega$ and all finite $\Gamma \sub A^\alpha_\beta$.
Furthermore, if $\U$ is a Ramsey ultrafilter, then $\FFF$ is called
{\it $\U-\alpha$--nice} (or simply {\it $\U$--nice})
iff given $f\in\omom$ one--to--one, there
is a club $D = D(f) \sub \omega_1$ such that for all $\beta \in D$
there exists $U\in\U$ with $f [U] \cap F_\beta^\gamma$ being finite
for all $\gamma \in A^\alpha_\beta$.
There is a two--way
interplay between niceness and $\U$--niceness (see 7.3 and 7.4): given 
$\FFF$ nice, we can,
in certain circumstances, construct $\U$ such that $\FFF$ is $\U$--nice;
on the other hand, after forcing with $\MM_\U$ where $\FFF$ is $\U$--nice,
$\FFF$ is still nice in the generic extension. This is the core of our
arguments, and guarantees the preservation of niceness in finite support
iterations with forcings of the form $\MM_\U$ as a factor.
\bigskip

{\capit Lemma 7.2.} (CH) {\it If $\FFF = \{ \F^\gamma = \{ F^\gamma_\beta
; \; \beta < \omega_1 \} ; \; \gamma < \alpha \}$ is $\alpha$--nice
and $\U$ is an ultrafilter,
then there is $\F^\alpha = \{ F^\alpha_\beta ; \; \beta < \omega_1 \} \sub \U$
such that $\FFF \cup \{ \F^\alpha \}$ is $(\alpha+1)$--nice.}
\sm

{\capit Lemma 7.3.} (CH) {\it If $\FFF$ is
nice, then there is a Ramsey ultrafilter
$\U$ such that $\FFF$ is $\U$--nice.}
\sm

{\capit Lemma 7.4.} {\it Assume $\U$ is a Ramsey ultrafilter, and
$\FFF$ is  $\U$--nice. Then $\forces_{\MM_\U}
`` \FFF$ is nice".}
\sm

{\capit Lemma 7.5.} {\it Let $\la \PP_\alpha , \dot\QQ_\alpha ; \;
\alpha < \delta \ra$, $\delta$ a limit ordinal, be a finite support iteration
of $ccc$ p.o.'s, and let $\dot\F^\alpha = \{
\dot F^\alpha_\beta ; \; \beta < \omega_1 \}$ be $\PP_\alpha$--names
for filter bases on $\omega$
such that $\forces_{\PP_\alpha} `` \dot\FFF^{\alpha + 1}:
= \{ \dot \F^\gamma ; \; \gamma < \alpha + 1\}$ is $(\alpha + 1
)$--nice$"$ for $\alpha
<\delta$. Then $\forces_{\PP_\delta}
`` \dot\FFF^\delta = \bigcup_{\alpha < \delta } \dot \FFF^{\alpha +1}$
is $\delta$--nice$"$.}
\bigskip


{\it Proof of Lemma 7.2.}
Let $\{ X_\beta ; \; \beta < \omega_1 \}$ enumerate $\omoms$
so that each $X \in\omoms$ occurs uncountably often.
Also let $\{ f_\beta ; \; \beta < \omega_1 \}$ enumerate
the one--to--one functions of $\omom$.
For $X_\eta $ and $\la f_\zeta \ra_{\zeta < \eta}$, 
let $C_{\eta} = C(X_\eta
, \la f_\zeta \ra_{\zeta < \eta} )$ witness the niceness of
$\FFF$. Without loss $\min (C_{\eta } )> \eta $.
It suffices to construct sets $F_\beta^\alpha \in \U$ such that for
all $\eta  < \beta$ with $\beta \in C_{\eta }$ and all
finite $\Gamma \sub A_\beta^\alpha \cup \{ \alpha \}$ and 
finite $Z \sub \eta$, we have
$| X_\eta \sem \bigcup_{\zeta \in Z} \bigcup_{\gamma \in \Gamma}
f^{-1}_\zeta (F^\gamma_\beta ) | = \omega$.
$F^\alpha_0$ is any member of $\U$.

Assume the $F_{\beta'}^\alpha$ have been constructed for $\beta '< \beta$.
Let $\{ \eta_k ; \; k\in\omega\}$ enumerate the $\eta < \beta$ 
with $\beta \in C_{\eta }$. Also, let $\{ \zeta_i ; \; i\in\omega\}$
enumerate $\beta$. By the properties of the $C_{\eta}$ with
$\beta \in C_{\eta}$ and standard thinning arguments, we can find
sets $X_k ' \sub X_{\eta_k}$ such that
\sm

\item{(i)} for all $\zeta \in \eta_k$ and $\gamma \in A_\beta^\alpha$ we have
that $X_k ' \cap f_\zeta^{-1} (F_\beta^\gamma)$ is finite;

\item{(ii)} the $X_k '$ are pairwise disjoint;

\item{(iii)} given $k$ and $i$, either $f_{\zeta_i}$ is almost equal to
some $f_{\zeta_j}$ with $j < i$ on the set $X_k '$ or $f_{\zeta_i}
[X_k ']$ is almost disjoint from $f_{\zeta_j} [\bigcup_\ell X_\ell ']$ for
all $j < i$.
\sm

\no Now choose infinite sets
$X^0_k$ and $X^1_k$ such that $X^0_k \cup X^1_k = X_k '$
and $X^0_k \cap X^1_k = \em$. Put $X^0 = \bigcup_k X^0_k$ and
$X^1 = \bigcup_k X^1_k$. Build disjoint sets $Y^0 = \bigcup_k Y^0_k
$ and $Y^1 = \bigcup_k Y^1_k$ as follows.
$Y^0_0 = f_{\zeta_0} [X^0] , Y^1_0 = f_{\zeta_0} [X^1] , ... ,
Y^0_k = f_{\zeta_k} [X^0] \sem \bigcup_{ j < k} Y^1_j,
Y^1_k = f_{\zeta_k} [X^1] \sem \bigcup_{ j < k} Y^0_j, ...$
We have either $Y^0 \not\in\U$ or $Y^1 \not\in\U$; without loss
the former holds, and we let $F=F_\beta^\alpha = \omega \sem Y^0 \in \U$.
$F$ is as required because we can now show by
induction on $j$  that $f_{\zeta_j}^{-1} (F) \cap X^0_{
k}$ is finite 
for all $j$ and $k$. 
This completes the proof of 7.2. $\qed$
\sm


{\it Proof of Lemma 7.3.} 
Let $\{ f_\beta ; \; \beta < \omega_1 \}$ enumerate
the one--to--one functions in $\omom$, and build a tower
$\{ U_\beta ; \; \beta < \omega_1 \}$ which generates a Ramsey
ultrafilter $\U$. Guarantee that $\U$ will be Ramsey in the successor
steps of the construction. If $\beta$ is a limit such that
$\beta \in \bigcap_{\theta < \beta} C( U_\theta , \la f_\zeta
\ra_{\zeta < \theta} )$, we can choose $U_\beta$ such that $U_\beta
\sub^* U_\theta$ for all $\theta <\beta$ and such that $U_\beta \cap
f_\theta^{-1} (F_\beta^\gamma )$ is finite for all $\theta < \beta$
and $\gamma \in A^\alpha_\beta$; otherwise let $U_\beta$ be any set almost
included in all $U_\theta$'s where $\theta < \beta$. Thus, if $f = f_\eta$,
$D(f) =  \{ \beta > \eta ; \; \beta \in \bigcap_{\theta < \beta}
C(U_\theta , \la f_\zeta \ra_{\zeta < \theta} ) \}$ is a diagonal
intersection of clubs, and thus a club, and witnesses $\U$--niceness. 
$\qed$
\sm


{\it Proof of Lemma 7.4.}
Let $\dot X$
and $\dot f^j$, $j\in\omega$, be $\MM_\U$--names for objects in $\omoms$ 
and for one--to--one functions in
$\omom$, respectively. Let $\dot x$ be the name for the increasing
enumeration of $\dot X$. Replacing $\dot X$ by a name for a subset of
$\dot X$, if necessary, we may assume that
\sm

\ce{$\forces_{\MM_\U} `` \dot f^j ( \dot x
(n)) \geq \dot m (n) "$ for all $j$ and all $n\geq j$, 
\hskip 1truecm $(\star)$}
\sm

\no where $\dot m$ denotes the canonical name for the Mathias--generic.
Since $\MM_\U$ and $\LL_\U$ are forcing equivalent
(see $\S$ 3), we can think of $\MM_\U$ as forcing with Laver trees $T$
such that the set of successors lies in $\U$ for every node above
the stem. Fix $n\in\omega$ and $j \leq n$.
Set
$$A(n,j) = \{ \sigma \in \omlom ; \; 
\hbox{ some } T^j_{\sigma,n} =
T \in \LL_\U \hbox{ with } stem(T) = \sigma
\hbox{ decides the value } 
\dot f^j( \dot x (n) )  \}.$$
Furthermore put
$$B(n,j) =\{ \sigma \in \omlom  ;\; \sigma\not\in A(n,j) \hbox{ and 
 } U^j_{\sigma , n}:= \{ k \in \omega ; \; \sigma \ha\la k\ra\in A(n,j) \}
\in \U  \}.$$
Call a triple $(\sigma,n,j)$ {\it relevant} iff $\sigma
\in B(n,j)$. For relevant triples $(\sigma , n,j)$ define
$f = f_{\sigma ,n}^j : U^j_{\sigma , n} \to \omega$ by
\sm

\ce{$f (k) =$ the value forced to $\dot f^j(\dot x (n))$ by 
$T_{\sigma\ha\la k\ra , n}^j
 .$}
\sm

\no  Using that $\U$
is Ramsey and that $f$ cannot be constant on a set from $\U$
by the definition of $B(n , j)$, we may assume that
$f$ is one--to--one on $U^j_{\sigma , n}$, by pruning that
set if necessary.  By $\U$--niceness find a club $D(f)$ such that for all
$\beta \in D$ there is $U \in \U$ with $f[U] \cap F^\gamma_\beta$
being finite for all $\gamma \in A^\alpha_\beta$. Let $C$ be the
intersection of all $D(f^j_{\sigma , n})$ where $(\sigma , n , j)$
is relevant. We claim that 
\sm

\ce{$\forces_{\LL_\U} `` | \dot f^j [\dot X] \cap F^\gamma_\beta |
< \omega "$ for all $j\in\omega, \beta \in C$ and $\gamma
\in A^\alpha_\beta$. \hskip 1truecm $(\star\star)$}
\sm

\no Clearly this suffices to complete the proof of the Lemma.

To see $(\star\star)$, fix $j\in\omega, \beta \in C , \gamma \in
A^\alpha_\beta$ and $T\in \LL_\U$. Put $\ell : = \max \{ j , 
| stem (T) | \}$. We shall recursively construct $T' \leq T$ such that
\sm

\ce{$T' \forces_{\LL_\U} `` \dot f^j (\dot x (n)) \notin F_\beta^\gamma
"$ for all $n\geq\ell$. \hskip 1truecm $(\star\star\star)$}
\sm

\no Along the construction we shall guarantee that if $\sigma \in T'
\cap A(n , j)$ for some $n \geq \ell$, then
$T_\sigma ' : = \{ \tau \in T' ; \; \tau \sub \sigma$ or
$\sigma \sub \tau \} \leq T^j_{\sigma , n }$ and that the value
forced to $\dot f^j (\dot x (n))$ by $T^j_{\sigma , n}$
does not belong to $F^\gamma_\beta$. By $(\star)$ we see that
$stem(T) \notin A(n,j)$ for all $n\geq \ell$. Hence we can put
$stem(T)$ into $T'$. To do the recursion step, assume we put $\sigma
\supseteq stem (T)$ into $T'$. Again by $(\star)$, the set $N$
of all $n$ such that $\sigma \in A(n,j)$ is finite. The same holds
for the set $M$ of all $n$ such that $\sigma \in B(n,j)$. By definition
of $C$, we can find $U \in \U$ such that $U \sub U^j_{\sigma , n}$ and
$f^j_{\sigma , n} [U] \cap F^\gamma_\beta = \em$ for all
$n\in M$. Now put $\sigma\ha\la k \ra$ into $T'$ iff $k \in U$ and
$\sigma\ha\la k \ra \in T^j_{\sigma , n}$ for all $n \in N$ and
$\sigma\ha\la k \ra \notin A(n,j)$ for all $n \in \omega
\sem (N \cup M)$. Using again $(\star)$, it is easily seen that the set of
all $k$ satisfying these three clauses belongs to $\U$. This completes
the recursive construction of $T'$. It is now easy to see
that $T'$ indeed satisfies $(\star\star\star)$. $\qed$

\sm


{\it Proof of Lemma 7.5.}
Let $\dot X$ and $\dot f^j$, $j\in\omega$, be $\PP_\delta$--names for objects
in $\omoms$ and one--to--one functions in
$\omom$, respectively. First assume that $cf (\delta) = \omega$,
and that $\delta = \bigcup_n \delta_n$ (where the $\delta_n$ form the
sequence fixed before the definition of the $\bar A^\alpha$). Now 
construct $\PP_{\delta_n}$--names $\dot X_n$ and $\dot f^j_n$
which can be thought of as approximations to our objects
as follows. Step into $V_{\delta_n}$. Find a decreasing sequence
of conditions $\la p_{n,m} ; \; m\in\omega \ra\in \PP_{[\delta_n , \delta)}
$ such that $p_{n,m}$ decides the $m$--th element of $\dot X$
as well as $\dot f^j (k)$ for $j,k \leq m$. Let 
$X_n$ be the set of elements forced into $\dot X$ by this
sequence, and let $f_n^j$ be the function whose values are
forced to $\dot f^j$ by this sequence. The niceness of
the $\FFF^{\delta_n + 1} $ in the models $V_{\delta_n}$
provides us with  clubs $ C_n= C ( X_n ,
\la f^j_n \ra_j )$ as witnesses.

Back in $V$, we have
names $\dot C_n$ for these witnesses.
By $ccc$--ness find a club $C$ in the ground model
which is forced to be contained in all $\dot C_n$ by the trivial condition
of $\PP_\delta$. We claim that $C$ witnesses the niceness
of $\FFF^\delta$ in $V_\delta$. To see this,
take $\beta \in C, k \in \omega$
and $\Gamma \sub A^\delta_\beta$ finite. Also fix
$\ell\in\omega$ and $p\in\PP_{\delta}$. Find $n$
such that $p\in \PP_{\delta_n}$ and $\Gamma \sub \delta_n$.
Step into $V_{\delta_n} = V[G_{\delta_n}]$
where $p \in G_{\delta_n}$. Since $\Gamma \sub A^{\delta_n}_\beta$
by construction of the $A^\delta_\beta$, we
know that $|X_n  \sem \bigcup_{j<k} \bigcup_{\gamma \in \Gamma} ( f_n^j)^{-1}
( F^\gamma_\beta) | = \omega$. Hence we can find $i > \ell$
in this set and $m$ large enough so that
\sm

\ce{$p_{n,m} \forces_{\PP_{[\delta_n,\delta)}} `` i \in \dot X \sem \bigcup_{j<k} 
\bigcup_{\gamma \in \Gamma} (\dot f^j)^{-1}
( F^\gamma_\beta) "$.}
\sm

\no Thus we see that
\sm

\ce{$\forces_{\PP_\delta}
`` | \dot X \sem \bigcup_{j<k} \bigcup_{\gamma \in \Gamma} (\dot f^j)^{-1}
(\dot F^\gamma_\beta) | = \omega"$,}
\sm

\no as required.

Next assume that $cf (\delta) = \omega_1$, and that $\delta = \bigcup_\zeta
\delta_\zeta$. Find $\alpha < \delta$ such that $\dot X$ and
$\dot f^j$ are $\PP_\alpha$--names, step into $V_\alpha$, and
let $X = \dot X [G_\alpha]$, $f^j = \dot f^j [G_\alpha]$. 
By niceness of the $\FFF^{\delta_\zeta + 1}
$ for $\delta_\zeta > \alpha$, we get
$\PP_{\delta_\zeta}$--names for clubs, $\dot C_\zeta = \dot C
(X , \la f^j\ra_j )$. Without loss, we can assume that
$C_\zeta = \dot C_\zeta [G_{\delta_\zeta}]
\in V_\alpha$ (by $ccc$--ness). Then $C = \{ \beta < \omega_1 ; \;
\beta \in \bigcap_{\zeta < \beta} C_\zeta \}$ is a diagonal intersection
of clubs, and thus club,
and is easily seen to witness the niceness of 
$\dot\FFF^\kappa [G_\kappa]$ (use again the definition of $\bar A^\delta$).
This completes the proof of the Lemma. $\qed$

\bigskip


{\it Proof of Theorem 8.} Let $E = \{ \alpha < \omega_2 ; \;
cf(\alpha) = \omega_1 \}$. We shall start with a model $V$ which satisfies
$GCH$ and additionally $\diamondsuit_{\omega_2} (E)$. The latter
is used as a guessing principle to
ensure that we took care of every ultrafilter along the iteration. For example
we could take $V = L$ (see [De, chapter IV, Theorem 2.2]).
Then we perform a finite support iteration $\la \PP_\alpha , \dot\QQ_\alpha ;
\; \alpha < \omega_2 \ra$ of $ccc$ p.o.'s over $V$.
We think of the diamond sequence as acting on the product $\omega_2
\times \PP_{\omega_2}$; more explicitly, we use a sequence
$\la S_\alpha \sub\alpha \times \PP_\alpha ; \; \alpha \in E \ra$ such that
for all $T \sub \omega_2 \times \PP_{\omega_2}$, the set $\{ \alpha \in E
; \; T \cap (\alpha\times \PP_{\alpha})
= S_\alpha \}$ is stationary. This we can do since the initial segments
$\PP_\alpha$ of the iteration will have size $\omega_1$.
Furthermore we shall have a $\PP_{\omega_2}$--name $\dot f$ for a bijection
between $\omega_2$ and $\omoms$ such that for all $\alpha \in E$,
we have
\sm

\ce{$\forces_{\PP_\alpha} `` \dot f \re \alpha$ is a bijection
between $\alpha$ and $\omoms \cap V[\dot G_\alpha] "$.}
\sm

\no The existence of such a name is, again, straightforward.

The details of the construction are as follows.
In $V_\alpha$, we shall have \sm

\item{(a)} a Ramsey ultrafilter $\U_\alpha$ such that $\dot\QQ_\alpha
[G_\alpha] = \MM_{\U_\alpha}$; 

\item{(b)} a filter base $\F^\alpha = \{ F^\alpha_\beta ; \; \beta
< \omega_1 \}$ such that $\FFF^{\alpha + 1} := \{ \F^\gamma ; \; \gamma
\leq \alpha \}$ is both nice and $\U_\alpha$--nice.
\sm

\no Let $\alpha $ be arbitrary. By either Lemma 7.4 or 7.5 and induction,
$\FFF^\alpha = \{ \F^\gamma ; \; \gamma < \alpha \}$ is nice 
in $V_\alpha$. In case $S_\alpha$ is a $\PP_\alpha$--name for a subset
of $\alpha$ and
\sm

\ce{$\forces_{\PP_\alpha} ``  \dot f [S_\alpha]$ is an ultrafilter",}
\sm

\no we let $\V = \dot f [S_\alpha] [G_\alpha]\in V_\alpha$; 
otherwise $\V$ is an
arbitrary ultrafilter of $V_\alpha$. By Lemma 7.2 find $\F^\alpha
\sub \V$ such that $\FFF^{\alpha + 1}$ is nice. 
Then apply Lemma 7.3 to get $\U_\alpha$
such that $\FFF^{\alpha + 1}$ is $\U_\alpha$--nice. 
This completes the construction.

It remains to see that $V_{\omega_2}$ is as required. $\cc =\ss = \omega_2$
is immediate because all the factors of the iteration are of the
form $\MM_\U$ for some Ramsey ultrafilter. To see $\pi\pp (\V)
= \omega_1$ for every ultrafilter $\V$, take a
$\PP_{\omega_2}$--name $T \sub \omega_2 \times \PP_{\omega_2}$
such that $\dot f [T] [G_{\omega_2}] = \V$; without loss
$\forces_{\PP_{\omega_2}} `` \dot f [T]$ is an ultrafilter".
We easily get a club $C \sub\omega_2$ such that for all $\alpha
\in C \cap E$, we have that $T\cap (\alpha \times \PP_\alpha)$ is a name
and
\sm

\ce{$\forces_{\PP_\alpha} `` 
\dot f [T \cap (\alpha\times \PP_\alpha)]$ is an ultrafilter in
$V[\dot G_\alpha]$".}
\sm

\no Hence we find $\alpha \in C \cap E$ with $T \cap (\alpha \times
\PP_\alpha) = S_\alpha$. This means that in $V_\alpha$, we chose
$\F^\alpha \sub \dot f [S_\alpha] [G_\alpha]$ such that $\F^\alpha$
had no pseudointersection in $V_{\omega_2}$. Since
$\F^\alpha \sub \dot f [T] [G_{\omega_2}] =\V$ has size $\omega_1$, $\pi\pp (\V)
= \omega_1$ follows. $\qed$
\bigskip

{\capit Remark 7.6.} If one cares only about $P$--points $\U$,
then the conclusion of Theorem 8 is much easier to prove because {\it
niceness} can be replaced by a simpler notion. Also, A. Dow
has remarked that $ \ss = \omega_2$ and $\pi\pp (\U) = \omega_1$
for all $P$--points $\U$ is true in Dordal's factored Mathias
real model [Do], and the referee has pointed out that one of the models
of [BlS 1] even satisfies $\ss = \omega_2$ and $\chi (\U) =
\omega_1$ for all $P$--points $\U$. The latter is so
because the forcing construction increases $\ss$ and is
$P$--point--preserving [BlS 1, Theorems 3.3 and 5.2]. The former
holds because in Dordal's model all $\omega_1$--towers
are preserved along the iteration (for the successor step,
one uses a result of Baumgartner [Do, Theorem 2.2], saying that Mathias
forcing does not destroy any towers; the limit step is taken care
of by the type of iteration used [Do, Lemma 4.2]). An easy reflection
argument shows each $P$--point $\U$ contains such an $\omega_1$--tower,
and $\pi\pp (\U) = \omega_1$ follows. We do not know whether Dordal's
model even satisfies $\pi\pp (\U) = \omega_1$ for all ultrafilters
$\U$. For this, one would have to extend Baumgartner's
result quoted above to filter bases. However, our construction
is more general, for slight modifications in the proof show the
consistency of the statement in Theorem 8 with large continuum;
more explicitly:
\bigskip

{\capit Remark 7.7.} Let $\kappa \geq \cc$ be a regular cardinal in $V
\models \diamondsuit_{\kappa^+} (E)$.
Then there is a generic extension of $V$ which satisfies $\pi\pp (\U)
= \kappa$ for all ultrafilters $\U$ and $\ss = \cc = \kappa^+$.
To see this, simply replace $\omega$, $\omega_1$ and $\omega_2$
by $< \kappa$, $\kappa$ and $\kappa^+$ (respectively) in the
above proof, and change the definitions of $\bar A^\alpha$ and
niceness accordingly. Then Lemmata 7.4 and 7.5 still hold 
(with a modified proof, of course) and Lemmata  7.2 and 7.3
are true if the assumption is changed to $MA + \cc = \kappa$.
This means that along the iteration we also have to force $MA$
cofinally often with $ccc$ p.o.'s of size $< \kappa$. This is no
problem since it can be shown (with an argument similar to the modified
proof of Lemma 7.5) that such p.o.'s preserve (the modified) niceness.
We leave details to the reader.

\Bigskip

\vfill\eject


{\dunhg 8. Questions with comments}
\Smallskip

\no There are numerous interesting questions connected with the cardinals
we have studied which are still open.
\bigskip

\item{(1)} {\it  Does $\chi_\sigma (\U) = \chi (\U)$ for all 
ultrafilters $\U$?}
\sm

\no We note that $\chi (\U) = \chi_\sigma (\U)$ as long as $\chi (\U)
< \omega_\omega$; furthermore, 
$\chi (\U) = \chi_\sigma (\U)$ in the absence of $0^\sharp$ (these
remarks are due to W. Just). 
\bigskip

\item{(2)} (Vojt\'a\v s, cf. [Va, Problem 1.4]) 
{\it Does $\rr = \rr_\sigma$?}
\sm

\no This problem is connected with Miller's question
whether $cf(\rr) = \omega$ is consistent (see [Mi, p. 502] and [Mi 1,
Problem 3.4]).
\bigskip

\item{(3)} {\it Does $\rr_\sigma = \min_\U \pi\chi_\sigma (\U)$?}
\sm

\no A negative answer would provide us with a dual form of Theorem 8,
and rescue some of the symmetry lost in $\S$ 7.
\bigskip

\item{(4)} {\it Can $\pi\pp(\U)$ be consistently singular?}
\sm

\no Let us recall ($\S$ 1) that $\pp (\U)$ is regular and 
notice that $\pi\chi (\U)$ and $
\chi (\U)$ are consistently singular --- simply add $\omega_{\omega_1}$
Cohen reals or see $\S\S$ 5 and 6! So $\pi\pp (\U)$ is the only cardinal for which this
question is of interest. 
Furthermore, we may ask whether $cf (\pi\pp (\U)) \geq \pp (\U)$.
(Note this is true for $\pi\chi (\U)$ and $\chi (\U)$, see
$\S$ 1.)
The only information we have about $cf (\pi\pp (\U))$ is given
in 1.7 and 1.8.
The problem seems connected with Vaughan's problem
concerning the possible singularity of $\ss$ (cf. [Va, Problem 1.2]).
\bigskip

\item{(5)} (Spectral problem at regulars) {\it Assume $\cc = \omega_3$
and there is an ultrafilter $\U$ with $\chi (\U) = \omega_1$.
Does this imply there is an ultrafilter $\V$ with $\chi (\V)
= \omega_2$? With $\pi\chi (\V)=\omega_2$?}
\sm

\no (Of course, this is just the smallest interesting case of a much more
general problem.) Note that the assumptions imply that there are ultrafilters
$\U$ and $\V$ with $\pi\chi (\U) = \chi (\U) = \omega_1$ and 
$\pi\chi (\V) = \chi (\V) = \omega_3$ (Bell--Kunen [BK], see
also [vM, Theorem 4.4.3]). By Theorem 6, we know there is not necessarily
an ultrafilter $\W$ with $\pi\chi (\W) = \chi (\W) = \omega_2$.
Of course, there is a corresponding problem on $\pp$ and $\pi\pp$.
Finally, similar questions can be asked about special classes of
ultrafilters. For example, it would be interesting to know
what can be said about the spectrum of possible characters of $P$--points.
\bigskip

\item{(6)} (Spectral problem at singulars) {\it Let $\lambda$ be
singular (of uncountable cofinality). Assume that \Spec{\chi} is
cofinal in $\lambda$. Does $\lambda \in$ \Spec\chi ? Similar
question for \Spec{\pi\chi}. What about $\lambda^+$?}
\sm

\no The only (partial) results we have in this direction are
Proposition 5.1 (a), Theorem 4 (c) and Proposition 6.3.
\bigskip

\item{(7)} {\it Let $R$ be a set of cardinals of uncountable cofinality in
$V \models GCH$. Show there is a generic extension of $V$ which has 
ultrafilters $\U$ with $\chi (\U) = \pi\chi (\U) = \lambda$ for 
each $\lambda \in R$.}
\sm

\no For regulars, this was done in Theorem 5. For singulars it was done
{\it separately} for $\pi\chi$ and $\chi$ in Corollaries 5.5 and
6.1. We don't know how to do it simultaneously. Note, however,
that given a single singular cardinal $\lambda$ of uncountable
cofinality in $V \models CH$, we can always force an ultrafilter with
$\chi (\U) = \pi\chi (\U) = \lambda$: simply add $\lambda$ Cohen reals;
then, in fact, all ultrafilters $\U$ satisfy $\chi (\U)
= \pi\chi (\U) = \lambda$.
\bigskip

\item{(8)} {\it Is there, in $ZFC$, an ultrafilter $\U$ with
$\pi\chi (\U) = \chi (\U)$?}
\sm

\no By the result of Bell and Kunen ([BK], [vM, Theorem 4.4.3]),
this is true if $\cc$ is regular. The Bell--Kunen model [BK]
which has no ultrafilter $\U$ with $\pi\chi (\U) = \cc$ has one
with $\pi\chi (\U) = \chi (\U) = \omega_1$ instead. The dual
question, whether there is an ultrafilter $\U$ with $\pp (\U) =
\pi\pp (\U)$, is independent, by the second author's $P$--point
independence Theorem [Sh]. However, we may still ask whether one
always has an ultrafilter $\U$ with $\pp ' (\U) = \pi\pp (\U)$.
\Bigskip

\vfill\eject


{\dunhg References}
\Smallskip

\itemitem{[BS]} {\capit B. Balcar and P. Simon,} {\it On minimal
$\pi$--character of points in extremally disconnected compact
spaces,} Topology and its Applications, vol. 41 (1991), pp.
133-145.
\sm
\itemitem{[Ba]} {\capit T. Bartoszy\'nski,} {\it On covering
of real line by null sets,} Pacific Journal of Mathematics,
vol. 131 (1988), pp. 1-12.
\smallskip
\itemitem{[BJ]} {\capit T. Bartoszy\'nski and H. Judah,}
{\it Measure and category --- filters on $\omega$,} in: Set Theory
of the Continuum (H. Judah, W. Just and H. Woodin, eds.),
MSRI Publications vol. 26 (1992), Springer, New York, pp. 175-201.
\smallskip
\itemitem{[BJ 1]} {\capit T. Bartoszy\'nski and H. Judah,}
{\it Set theory --- On the structure of the real line,} A K Peters,
Wellesley, 1995.
\sm
\itemitem{[B]} {\capit J. Baumgartner,} {\it Iterated forcing,}
Surveys in set theory (edited by A.R.D. Mathias), Cambridge
University Press, Cambridge, 1983, pp. 1-59.
\smallskip
\itemitem{[Be]} {\capit M. G. Bell,} {\it On the combinatorial
principle $P(c)$,} Fundamenta Mathematicae, vol. 114 (1981),
pp. 149-157.
\smallskip
\itemitem{[BK]} {\capit M. Bell and K. Kunen,} {\it On the
$\pi$--character of ultrafilters,} C. R. Math. Rep.
Acad. Sci. Canada, vol. 3 (1981), pp. 351-356.
\sm
\itemitem{[Bl]} {\capit A. Blass,} {\it Near coherence of filters,
I: cofinal equivalence of models of arithmetic,}
Notre Dame Journal of Formal Logic, vol. 27 (1986), pp.
579-591.
\sm
\itemitem{[Bl 1]} {\capit A. Blass,} {\it Selective ultrafilters 
and homogeneity,} Annals of Pure and Applied Logic,
vol. 38 (1988), pp. 215-255.
\smallskip
\itemitem{[Bl 2]} {\capit A. Blass,} {\it Simple cardinal
characteristics of the continuum,} in: Set Theory of the Reals
(H. Judah, ed.), Israel Mathematical Conference Proceedings, vol. 6,
1993, pp. 63-90.
\smallskip
\itemitem{[BlM]} {\capit A. Blass and H. Mildenberger,}
{\it On the cofinality of ultraproducts,} preprint.
\sm
\itemitem{[BlS]} {\capit A. Blass and S. Shelah,} {\it Ultrafilters with
small generating sets,} Israel Journal of Mathematics, vol. 65 (1989),
pp. 259-271.
\smallskip
\itemitem{[BlS 1]} {\capit A. Blass and S. Shelah,}
{\it There may be simple $P_{\aleph_1}$-- and $P_{\aleph_2}$--points
and the Rudin--Keisler ordering may be downward directed,}
Annals of Pure and Applied Logic, vol. 33 (1987),
pp. 213-243.
\smallskip
\itemitem{[Br]} {\capit J. Brendle,} {\it Strolling through
paradise,} Fundamenta Mathematicae, vol. 148 (1995), pp. 1-25.
\smallskip
\itemitem{[Ca]} {\capit R. M. Canjar,} {\it Countable
ultrapoducts without $CH$,} Annals of Pure and Applied Logic,
vol. 37 (1988), pp. 1-79.
\sm
\itemitem{[Ca 1]} {\capit R. M. Canjar,} {\it Mathias forcing which does
not add dominating reals,} Proceedings of the American Mathematical
Society, vol. 104 (1988), pp. 1239-1248.
\smallskip
\itemitem{[Ca 2]} {\capit R. M. Canjar,} {\it On the generic existence
of special ultrafilters,} Proceedings of the American Mathematical
Society, vol. 110 (1990), pp. 233-241.
\smallskip
\itemitem{[De]} {\capit K. Devlin,} {\it Constructibility,}
Springer, Berlin, 1984.
\smallskip
\itemitem{[Do]} {\capit P. L. Dordal,} {\it A model in which
the base--matrix tree cannot have cofinal branches,}
Journal of Symbolic Logic, vol. 52 (1987), pp. 651-664.
\sm
\itemitem{[Do 1]} {\capit P. L. Dordal,} {\it Towers in
$\omoms$ and $\omom$,} Annals of Pure and Applied Logic, vol. 45
(1989), pp. 247-276.
\sm
\itemitem{[El]} {\capit E. Ellentuck,} {\it A new proof that analytic sets
are Ramsey,} Journal of Symbolic Logic, vol. 39 (1974), pp. 163-165.
\sm
\itemitem{[GRSS]} {\capit M. Goldstern, M. Repick\'y, S. Shelah
and O. Spinas,} {\it On tree ideals,} Proceedings of the American
Mathematical Society, vol. 123 (1995), pp. 1573-1581.
\smallskip
\itemitem{[GS]} {\capit M. Goldstern and S. Shelah,} {\it Ramsey 
ultrafilters and the reaping number,} Annals of Pure and Applied
Logic, vol. 49 (1990), pp. 121-142.
\smallskip
\itemitem{[Je]} {\capit T. Jech,} {\it Set theory,} Academic Press,
San Diego, 1978.
\smallskip
\itemitem{[Je 1]} {\capit T. Jech,} {\it Multiple forcing,}
Cambridge University Press, Cambridge, 1986.
\smallskip
\itemitem{[JS]} {\capit H. Judah and S. Shelah,} {\it 
$\Delta_2^1$--sets of reals,} Annals of Pure and Applied Logic,
vol. 42 (1989), pp. 207-223.
\smallskip
\itemitem{[Ku]} {\capit K. Kunen,} {\it Set theory,} North-Holland,
Amsterdam, 1980.
\sm
\itemitem{[LR]} {\capit G. Lab\c edzki and M. Repick\'y,} {\it
Hechler reals,} Journal of Symbolic Logic, vol. 60 (1995),
pp. 444-458.
\smallskip
\itemitem{[Lo]} {\capit A. Louveau,} {\it Une m\'ethode topologique
pour l'\'etude de la propri\'et\'e de Ramsey,} Israel Journal of
Mathematics, vol. 23 (1976), pp. 97-116.
\smallskip
\itemitem{[M]} {\capit P. Matet,} {\it Combinatorics and forcing
with distributive ideals,} Annals of Pure and Applied Logic,
vol. 86 (1997), pp. 137-201.
\sm
\itemitem{[Ma]} {\capit A. R. D. Mathias,} {\it Happy families,}
Annals of Mathematical Logic, vol. 12 (1977), pp. 59-111.
\smallskip
\itemitem{[Mi]} {\capit A. Miller,} {\it A characterization of the
least cardinal for which the Baire category theorem fails,} Proceedings
of the American Mathematical Society, vol. 86 (1982), pp. 498-502.
\smallskip
\itemitem{[Mi 1]} {\capit A. Miller,} {\it Arnie Miller's problem list,}
in: Set Theory of the Reals (H. Judah, ed.), Israel Mathematical
Conference Proceedings, vol. 6, 1993, pp. 645-654.
\smallskip
\itemitem{[Ny]} {\capit P. Nyikos,} {\it Special ultrafilters
and cofinal subsets of $\omom$,} preprint.
\sm
\itemitem{[Pl]} {\capit S. Plewik,} {\it On completely Ramsey sets,}
Fundamenta Mathematicae, vol. 127 (1986), pp. 127-132.
\smallskip
\itemitem{[Sh]} {\capit S. Shelah,} {\it Proper forcing,}
Lecture Notes in Mathematics, vol. 940, Springer--Verlag,
New York, 1982.
\sm
\itemitem{[SS]} {\capit S. Shelah and J. Stepr\=ans,}
{\it Maximal chains in $\omom$ and ultrapowers of the integers,}
Archive for Mathematical Logic, vol. 32 (1993), pp. 305-319, and
vol. 33 (1994), pp. 167-168.
\sm
\itemitem{[St]} {\capit J. Stepr\=ans,} {\it Combinatorial
consequences of adding Cohen reals,} in: Set Theory of the Reals
(H. Judah, ed.), Israel Mathematical Conference Proceedings, vol.
6, 1993, pp. 583-617.
\smallskip
\itemitem{[vD]} {\capit E. K. van Douwen,} {\it
The integers and topology,} Handbook of set--theoretic
topology, K. Kunen and J. E. Vaughan (editors),
North--Holland, Amsterdam, 1984, pp. 111-167.
\smallskip
\itemitem{[vM]} {\capit J. van Mill,} {\it An introduction to
$\beta\omega$,}  Handbook of set--theoretic
topology, K. Kunen and J. E. Vaughan (editors),
North--Holland, Amsterdam, 1984, pp. 503-567.
\sm
\itemitem{[Va]} {\capit J. Vaughan,} {\it Small
uncountable cardinals and topology,} in: Open problems
in topology (J. van Mill and G. Reed, eds.),
North--Holland, 1990, pp. 195-218.
\smallskip

\vfill\eject\end